\documentclass{article}
\usepackage{amssymb}
\newtheorem{theorem}{Theorem}
\newtheorem{definition}{Definition}
\newtheorem{proposition}{Proposition}
\newtheorem{lemma}{Lemma}

\title{Dirac Type Operators for Arithmetic Subgroups of
Generalized Modular Groups}
\author{E. Bulla \thanks{Department of Mathematics, Louisiana State University, Baton Rouge, LA 70803, USA. Financial support by the SILO/SURF grant from the Arkansas Department of Education gratefully acknowledged. E-mail:
{\tt ebulla
math.lsu.edu}}
\and
D.~Constales \thanks{Department of
Mathematical Analysis, Ghent University, Building
S-22, Galglaan 2, B-9000 Ghent, Belgium. Financial support from BOF/GOA 01GA0405 of Ghent University gratefully acknowledged.  E-mail:
{\tt dc
cage.UGent.be}}
\and R. S. Krau{\ss}har \thanks{Department of
Mathematical Analysis, Ghent University, Building
S-22, Galglaan 2, B-9000 Ghent, Belgium. E-mail:
{\tt krauss
cage.UGent.be}}
\and John Ryan \thanks{Department of Mathematics, University of Arkansas, Fayetteville, AR 72701, USA {\tt jryan
uark.edu}}}
\begin{document}
\maketitle
\date
\begin{abstract}
Fundamental solutions of Dirac type operators are introduced for a class of conformally flat manifolds. This class consists of manifolds obtained by factoring out the upper half-space of $\mathbb{R}^n$ by arithmetic subgroups of generalized modular groups. Basic properties of these fundamental solutions are presented together with associated Eisenstein series.  
\end{abstract}
\section{Introduction}

\ A natural generalization to $\mathbb{R}^{n}$ of the classical Cauchy-Riemann operator has proved to be the euclidean Dirac operator. Associated to this operator is a Cauchy Integral Formula and other natural analogues of basic results from one variable complex analysis. See for instance \cite{bds} and elsewhere. Further the euclidean Dirac operator has been used in understanding boundary value problems and aspects of classical harmonic analysis in $\mathbb{R}^{n}$. See for instance \cite{GS2,sijuewu} and elsewhere. This analysis together with its applications is known as Clifford analysis.

\ On the other hand Dirac operators have proved to be extremely useful tools in understanding geometry over spin manifolds. See for instance \cite{LM} and elsewhere. Basic aspects of Clifford analysis over spin manifolds have been developed in \cite{calderbank,cn,mitrea}. Further in \cite{KraRyan1,KraRyan2,KraRyan3,lr,ma,r85} and elsewhere it is illustrated that the context of conformally flat manifolds provide a useful setting for developing Clifford analysis.

\ Conformally flat manifolds are those manifolds which possess an atlas whose transition functions are M\"obius transformations. Under this viewpoint conformally flat manifolds can be regarded as higher dimensional generalizations of Riemann surfaces, as pointed out for example in \cite{ma,r85}. These types of manifolds have been studied extensively in a number of contexts. See for instance \cite{Kuiper,sy}.
Following the classical work of  N. H. Kuiper \cite{Kuiper}, one can construct
a whole family of examples of conformally flat manifolds by factoring out a
subdomain $U \subseteq \mathbb{R}^n$ by a Kleinian group $\Gamma$ acting totally
discontinuously on $U$.

\par\medskip\par

Simple examples of conformally flat manifolds include spheres, hyperbolas, real projective space,
cylinders, tori, and the Hopf manifolds $S^1 \times S^{n-1}$. In
\cite{lr,KraRyan1,KraRyan2} explicit Clifford analysis techniques have been developed for these manifolds.

\par\medskip\par

In this paper we treat some special examples of hyperbolic manifolds of higher
genus with spinor structure. The manifolds that we consider arise
from factoring out upper half-space in 
$\mathbb{R}^n$ by an arithmetic congruence group, $H$, of the generalized modular
group $\Gamma_p$. $\Gamma_p$ is the group that is generated
by $p$ translation matrices ($p < n$) and the inversion
matrix.  

The
associated manifolds are higher dimensional analogues
of those classical Riemann surfaces that arise from factoring out the complex
upper half-plane by the principal congruence subgroups of 
$SL(2,\mathbb{Z})$. In two real variables these 
are $k$-handled spheres. 

\par\medskip\par
 
Spinor sections on
these manifolds can be constructed from automorphic forms on $H$. Using these
automorphic forms we set up Cauchy kernels or fundamental solutions of the Dirac operators associated to these spin manifolds. 

\par\medskip\par

The basic theory of monogenic, Euclidean harmonic and more generally of
$k$-genic automorphic forms on this family of arithmetic groups is described
in \cite{KraHabil} using generalized Eisenstein series. 
However, as shown in \cite{KraHabil}, even in the monogenic case, the absolute convergence abscissa with respect to the generalized modular group of these generalized Eisenstein series is only $p < n-2$. Here we overcome this convergence problem for the cases $p=n-2$ and $p=n-1$ by adapting a classical trick of Hecke, mentioned for example in \cite{Freitag}.

\  Dirac operators associated to upper half space endowed with the hyperbolic metric and scalar perturbations of the hyperbolic metric and their associated hyperbolic Laplace operators have received steady attention. See for instance \cite{AL,a-ord,slberlin,sl,slnewnew,Hua,Leutw87,Leutwiler,qber}. In this paper we introduce generalized Eisenstein series, with respect to the groups considered here, that are solutions to the these hyperbolic differential operators. We then use these series to introduce fundamental solutions to a particular hyperbolic Dirac operator and also for the hyperbolic Laplacian with respect to the manifolds considered here. 
\par\medskip\par
The lay out of the paper is as follows.

In Section 2 we introduce the background on Clifford algebras and Clifford analysis that we shall need here. 

\ In Section 3 some of the geometry of these manifolds is described. In particular direct analogues of fundamental domains described for arithmetic groups in $2$ and $3$ real variables are introduced. The Ahlfors-Vahlen approach to describing M\"{o}bius transformations in $n$ real variables is used to introduce isometric spheres to describe fundamental domains and their associated conformally flat manifolds.

\ In Section 4 for the cases $p=1'\ldots, n-3$ generalized Eisenstein series are introduced over the universal
covering space, $H^+(\mathbb{R}^{n})=\{x\in \mathbb{R}^{n}:x_{n}>0\}$,
of these manifolds and the projection operator from
$H^+(\mathbb{R}^{n})$
to the manifold is used to produce non-trivial solutions to Dirac
type equations over a spinor bundle. In this section the Hecke trick is adapted to introduce
similar sections  for the cases $p=n-2$ and $p=n-1$.

In Section 5 we construct special variants of Poincar\'e type series which induce the explicit
Cauchy kernels for monogenic sections. Other generalized Poincar\'e series are used to explicitly determine fundamental solutions to higher order Dirac type operators on the manifolds. Basic properties of these fundamental solutions are investigated. In particular results mentioned in Section 3 are used to obtain Hardy space decompositions of the $L^{q}$ spaces of compact strongly Lipschitz hypersurfaces lying in the manifolds considered here. This is for $q\in (1,\infty)$. Further the techniques used here are adapted to introduce operators of Calderon-Zygmund type in this context.

\par\medskip\par

In Section 6 we develop the analogous results for $k$-hypergenic Eisenstein series and $k$-hyperbolic harmonic Eisenstein series.  This function class includes hypermonogenic Eisenstein series  when $k=n-2$ and hyperbolic harmonic Eisenstein series also when $k=n-2$. Analogous sections are set up over the corresponding conformally flat manifolds. In the second part of Section 6 fundamental solutions to the Dirac operator and hyperbolic Laplacian associated to these Eisenstein series are introduced and some of their basic properties are investigated. In particular we also provide a hypermonogenic Hardy space decomposition of the $L^{q}$ space of a strongly Lipschitz hypersurface of the manifolds considered here. Again $q\in (1,\infty)$. 

\par\medskip\par

{\bf Acknowledgement}. The authors are very thankful to Aloys Krieg for helpful discussions while preparing this paper. 

\section{Preliminaries}
\subsection{Clifford algebras}
In this subsection we will introduce the basic information on Clifford algebras that we need in this paper. We shall regard Euclidean space, $\mathbb{R}^{n}$, as being embedded in the real, $2^{n}$-dimensional, associative Clifford algebra, $Cl_{n}$, satisfying the relation $x^{2}=-\|x\|^{2}$ for each $x\in \mathbb{R}^{n}$. In terms of the standard orthonormal basis 
$e_1,\ldots,e_n$ of $\mathbb{R}^{n}$ this relation becomes the anti-commutation $e_i e_j + e_j e_i = - 2 \delta_{ij}$ and a basis for $Cl_{n}$ is given by 
\[1, e_{1},\ldots,e_{n},\ldots,e_{j_{1}}\ldots e_{j_{r}},\ldots, e_{1}\ldots e_{n}\]
where $j_{1}<\ldots j_{r}$ and $1\leq r\leq n$. Each non-zero vector in $\mathbb{R}^{n}$ has a multiplicative inverse, $x^{-1}=\frac{-x}{\|x\|^{2}}$. Up to a sign this is the Kelvin inverse of a non-zero vector. For $a=a_{0}+\ldots +a_{1\ldots n}e_{1}\ldots e_{n}\in Cl_{n}$ the nom of $a$ is defined to be $\|a\|=(a_{0}^{2}+\ldots a_{1\ldots n}^{2})^{\frac{1}{2}}$. The reversion anti-automorphism is defined by  $\sim:Cl_{n}\rightarrow Cl_{n}:\sim e_{j_{1}}\ldots e_{j_{r}}=e_{j_{r}}\ldots e_{j_{1}}$. For $a\in Cl_{n}$ we write $\tilde{a}$ for $\sim a$.

\ Any element $a \in Cl_n$ may be uniquely decomposed in the form $a = b + c e_n$, where $b,c \in Cl_{n-1}$. Based on this decomposition one defines the projection mappings $P: Cl_n \rightarrow Cl_{n-1}$ and $Q: Cl_n \rightarrow Cl_{n-1}$ by $Pa=b$ and $Qa=c$. Further we define $Q^{\star}a$ to be $e_{n}(Qa)e_{n}$. Note that if we define $\hat{a}$ to be $b-ce_{n}$ then 
\begin{equation}
P(a)=\frac{1}{2}(a+\hat{a})
\end{equation}
and
\begin{equation}
Q(a)=\frac{1}{2}(a-\hat{a}).
\end{equation}

\subsection{PDE's related to the Dirac operator in $\mathbb{R}^{n}$ and hyperbolic space}
{\bf Monogenic and $k$-genic functions.} 

\ The Dirac operator, $D$, in $\mathbb{R}^{n}$ is defined to be $ \frac{\partial }{\partial x_1} e_1 + \frac{\partial }{\partial x_2} e_2 + \cdots + \frac{\partial }{\partial x_n}e_n$. Suppose that $U$ is a domain in $\mathbb{R}^{n}$ and $f$ and $g$ are pointwise differentiable functions defined on $U$ and taking values in $Cl_{n}$. The function $f$ is called left monogenic or left Clifford holomorphic if it satisfies the equation $Df=0$ on $U$. Similarly $g$ is called right monogenic, or right Clifford holomorphic if it satisfies the equation $gD=0$ on $U$. Here $gD:=\Sigma_{j=1}^{n}\frac{\partial g(x)}{\partial x_{j}}e_{j}$.
 Due to the non-commutativity of $Cl_{n}$ for $n>1$, both
classes of functions do not coincide with each other. However, $f$ is left monogenic if and only if $\tilde{f}$ is right monogenic. The left and right fundamental solution to the $D$-operator is called the Euclidean Cauchy kernel and has the form $G_1(x-y) = \frac{1}{\omega_{n}}\frac{x-y}{\|x-y\|^n}$. Here $\omega_{n}$ is the surface area of the unit sphere in $\mathbb{R}^{n}$.
The Dirac operator factorizes the Euclidean Laplacian $\Delta = \sum_{j=1}^n \frac{\partial^2}{\partial x_j^2}$, viz $D^2 = - \Delta$. Every real component of a monogenic function is hence harmonic. More generally,
functions satisfying $D^k f = 0$ for a positive integer $k$ are called  left $k$-genic, while functions satisfying $gD^{k}=0$ are called right $k$-genic. Note that when $k$ is even, say $k=2l$ for some positive integer $l$, then $D^{2l}=(-\triangle)^{l}$ and in this case left and right $k$-genic functions coincide.

\ The fundamental solution to the operator $D^{k}$ for $k<n$ is $G_{k}(x-y)=C_{k}\frac{x-y}{\|x-y\|^{n-k+1}}$ when $k$ is odd and $G_{k}(x-y)=\frac{C_{k}}{\|x-y\|^{n-k}}$ when $k$ is even. Further $C_{k}$ is a real positive constant chosen so that $D^{k-1}G_{k}(x-y)=G_{1}(x-y)$. In what follows we shall write $G_{1}(x-y)$ simply as $G(x-y)$.

{\bf{Ahlfors-Vahlen matrices and iterated Dirac operators.}} 

\ Following for example \cite{a,EGM87}, M\"obius transformations in $\mathbb{R}^n$ can be
represented as
$$
m: \mathbb{R}^n \cup \{\infty\} \rightarrow \mathbb{R}^n \cup \{\infty\}:\;\; m(x)=(ax+b)(cx+d)^{-1}
$$
with coefficients $a,b,c,d$ from $Cl_n$ that can all be written as products of vectors from $\mathbb{R}^n$. Further 
$a\tilde{d}-b\tilde{c}\in \mathbb{R}\backslash\{0\}$ and $\tilde{a}c$, $\tilde{c}d$, $\tilde{d}b$, $\tilde{b}a \in \mathbb{R}^n$. These conditions are often called
Vahlen conditions.

\par\medskip\par

The set that consists of Clifford valued matrices $\left(\begin{array}{cc} a &
b \\ c & d \end{array}\right)$ whose coefficients satisfy the previously mentioned conditions is a group under matrix multiplication. It is called the general Ahlfors-Vahlen group. It is denoted by $GAV(\mathbb{R}^n)$. 

\par\medskip\par

The general linear Ahlfors-Vahlen group $GAV(\mathbb{R}^n)$ is a generalization of the general linear group $GL(2,\mathbb{R})$. The particular subgroup
$$
SAV(\mathbb{R}^n) =\{M \in GAV(\mathbb{R}^n)\;|\; a\tilde{d} - b \tilde{c} = 1\}
$$
is called the special Ahlfors-Vahlen group and it is purely generated by the inversion matrix and translation type matrices, as proved for instance in \cite{EGM87}. 

The projective special Ahlfors-Vahlen group is the group $$PSAV(\mathbb{R}^{n}) \stackrel{\sim}{=}SAV(\mathbb{R}^{n})/\{\pm I\}$$ where $I$ is the identity matrix.

The 
subgroup $SAV(\mathbb{R}^{n-1})$ of $SAV(\mathbb{R}^{n})$ has
the special property that it acts transitively on the upper half space $H^{+}(\mathbb{R}^{n})=\{x=x_{1}e_{1}+\ldots+x_{n}e_{n}\in\mathbb{R}^{n}:x_{n}>0\}$.

\par\medskip\par

Now assume that $m(x)$, $=M<x>$,  is a M\"obius transformation represented in the above mentioned form. It is shown in \cite{KraHabil} and elsewhere that if $f$ is a left $k$-genic  function in the  variable
$y=m(x)=(ax+b)(cx+d)^{-1}$, then the function $J_k(M,x)f(M<x>)$
is again left $k$-genic but now in the variable $x$.

Here
$$J_k(M,x)
= \left\{ \begin{array}{cc} \frac{\widetilde{cx+d}}{\|cx+d\|^{n-k+1}} &  \;\; k\;{\mbox odd}\\   \frac{1}{\|cx+d\|^{n-k}} & \;\;k\;{\mbox even}\end{array}\right.$$
In what follows we always restrict attention to the cases where $k < n$. This type of invariance for $k$-genic functions under M\"{o}bius transformations is seen as an automorphic invariance.

It may now be seen that if $f$ is a
function that is left $k$-genic on $H^{+}(\mathbb{R}^{n})$ then so is $J_k(M,x)f(M<x>)$, for any $M\in SAV(\mathbb{R}^{n-1})$. In what follows we shall write $J(M,x)$ for $J_{1}(M,x)$.

\par\medskip\par

{\bf $k$-hypergenic functions.}\\ 
Left $k$-hypergenic functions are defined as the null-solutions to the system
$$
M_{k}f:=D f + \frac{kQ^{\star} f}{x_n} = 0,\quad\quad x_n \neq 0,
$$
where $k \in \mathbb{R}$. This is a Hodge-Dirac equation for upper half space equipped with the metric $x_{n}^{\frac{-2k}{n-2}}
\sum\limits_{j=1}^{n} dx_{j}^{2}$. In the case $k=0$ these are precisely left monogenic functions. A similar definition can be given for right $k$-hypergenic functions.
\par\medskip\par

\ It is pointed out in \cite{qber} and elsewhere that
$$
-M_{k}^{2}f=(\triangle Pf-\frac{k}{x_{n}}\frac{\partial Pf}{\partial x_{n}})+(\triangle Qf-\frac{k}{x_{n}}\frac{\partial Qf}{\partial x_{n}}+\frac{k}{x_{n}^{2}}Qf).
$$

\ The operator $\triangle -\frac{k}{x_{n}}\frac{\partial}{\partial x_{n}}$ is the Laplacian for $H^{+}(\mathbb{R}^{n})$ with respect to the metric  $x_{n}^{\frac{-2k}{n-2}}\sum_{j=1}^{n}dx_{j}^{2}$. We will denote it by $\triangle_{k}$ and we shall call it the $k$-hyperbolic Laplacian. We will denote the operator $\triangle_{k}+\frac{k}{x_{n}^{2}}$ by $W_{k}$. Functions that are annihilated by the $k$-hyperbolic Laplacian will be called $k$-hyperbolic harmonic functions. When $k=n-2$ this operator is the hyperbolic Laplacian and functions that are annihilated by this operator are called hyperbolic harmonic functions. The operators $\triangle_{k}$ and $W_{k}$ are special cases of the Weinstein equation described in \cite{AL,Leutw87} and elsewhere. 

\ Under a M\"{o}bius transformation $y=M<x>=(ax+b)(cx+d)^{-1}$ a $k$-hypergenic function $f(y)$ is transformed to the $k$-hypergenic function
\begin{equation}
\label{kinv}
F(x):=K_{k}(M,x)f(M<x>),
\end{equation}

where $K_{k}(M,x)=\frac{\widetilde{cx+d}}{\|cx+d\|^{n-k}}$.
See for instance \cite{slnewnew,qber}.
The particular solutions associated to the case $k=n-2$ coincide with the null-solutions to the hyperbolic Hodge-Dirac operator with respect to the hyperbolic metric on upper half space. These are often called hyperbolic monogenic functions or simply hypermonogenic functions \cite{sl,Leutwiler}. Note that when $k$ is even then $K_{k}(M,x)=J_{k+1}(M,x)$.

\ As explained in \cite{sl,qber} the basic hypermonogenic kernels for $H^{+}(\mathbb{R}^{n})$ are given by
$$
p(x,y)=\frac{1}{\omega_{n}}x_{n}^{n-2}y_{n}^{n-1}\frac{(x-y)}{\|x-y\|^{n}}e_{n}\frac{(x-\hat{y})}{\|x-\hat{y}\|^{n}}
$$
and
$$
q(x,y)=y_{n}^{n-2}(\frac{1}{\|x-\hat{y}\|^{n-2}}G(x-y)+\frac{1}{\|x-y\|^{n-2}}G(x-\hat{y}))=y_{n}^{n-2}D_{x}H(x,y)
$$
where $D_{x}$ is the Dirac operator with respect to the variable $x$ and $H(x,y)=\frac{1}{(n-2)\omega_{n}}\frac{1}{\|x-y\|^{n-2}\|x-\hat{y}\|^{n-2}}$.

These kernels provide analogues of the Cauchy kernel $G(x-y)$ for hypermonogenic functions. More precisely if $f(x)$ is a left hypermonogenic function defined on a domain $U\subset{\cal{F}}_{p}[N]$ then
\begin{equation}
f(y)=P(\int_{\partial V}p(x,y)\frac{n(x)}{x_{n}^{n-2}}f(x)d\sigma(x))+Q(\int_{\partial V}q(x,y)n(x)f(x)d\sigma(x)
\end{equation}
where $V$ is a subdomain of $U$ whose closure is compact and lies in $U$. Further $y\in V$.

\ Using Equations (1) and (2) it is shown in \cite{sl} that Equation (4) can be rewritten as
$$
f(y)=y_{n}^{n-1}(\int_{\partial V}E(x,y)n(x)f(x)d\sigma(x)+\int_{\partial V}F(x,y)\hat{n}(x)\hat{f}(x)d\sigma(x))
$$
where $E(x,y)=\frac{2^{n-1}}{\|x-\hat{y}\|^{n-2}}G(x-y)$ and $F(x,y)=\frac{2^{n-1}}{\|x-y\|^{n-2}}G(\hat{x}-y)$.

\ It is shown in \cite{AL} that if $M\in SAV(\mathbb{R}^{n-1})$ then a $k$ hyperbolic harmonic function $f(y)$ is transformed to the $k$- hyperbolic harmonic function $L_{k}(M,x)f(M<x>)$ where $y=M<x>$ and $L_{k}(M,x)=\frac{1}{\|cx+d\|^{n-k-2}}$.

\ Following \cite{slberlin}  we also have by direct computation:

\begin{proposition}
 A function $f$ is left $k$-hypergenic if and only if $\frac{f}{x_{n}^{k}}e_{n}$ is left $-k$-hypergenic.
\end{proposition}

\ In \cite{AL} it is shown that if $u(y)$ is a solution to $W_{k}u=0$ then so is $L_{k}(M,x)u(M<x>)$ where $y=M<x>$. 

\ A function that is annihilated by the operator $W_{k}$ is an eigenfunction of the operator $x_{n}^{2}\triangle_{k}$ with eigenvalue $-k$.
 
\section{Subgroups of generalized modular groups, their fundamental domains and associated conformally flat manifolds}

\subsection{Arithmetic subgroups of the Ahlfors-Vahlen group}
Arithmetic subgroups of the special Ahlfors-Vahlen group that act totally discontinuously on $H^{+}(\mathbb{R}^{n})$ are for instance considered in \cite{EGM87,EGM90,EGM98} and for the three dimensional case in \cite{Maa49}. 

\ First let us introduce the ring
$$
{\cal{O}}_p := \sum_{A \subseteq  P(1,\ldots,p)} \mathbb{Z} e_A \quad \quad p \le  n-1.
$$
This ring of course lies in the subalgebra $Cl_{p}$. In what follows, let
$$
J:= \left(\begin{array}{cc} 0 & -1 \\ 1 & 0 \end{array}\right),T_{e_1} := \left(\begin{array}{cc} 1 & e_1 \\ 0 & 1  \end{array}\right),\ldots, T_{e_p}:= \left(\begin{array}{cc} 1 & e_p \\ 0 & 1  \end{array}\right).
$$
We recall, cf. \cite{KraHabil}:
\begin{definition} For $p < n$, the special hypercomplex modular groups are defined to be $\Gamma_p:= \langle J,T_{e_1},\ldots,T_{e_p}\rangle$. For a positive integer $N$ the associated principal
congruence subgroups of $\Gamma_p$ of level $N$\label{congsub} are then given by
$$
\Gamma_p[N] := \Big\{ \left(\begin{array}{cc} a & b\\ c & d \end{array}\right) \in \Gamma_p\;\Big|\; a-1,b,c,d-1 \in N  {\cal{O}}_p\Big\}.
$$
\end{definition}

\ As the group $SL(2,\mathbb{Z})$ is generated by $\left(\begin{array}{cc} 0 & -1 \\ 1 & 0\end{array}\right)$ and $\left(\begin{array}{cc} 1 & 1\\ 0 & 1\end{array}\right)$ it can be seen that $\Gamma_{p}$ is a natural generalization to $\mathbb{R}^{n}$ of $SL(2,\mathbb{Z})$. It follows that the group $\Gamma_{p}[N]$ is a natural generalization of arithmetic subgroups of $SL(2,\mathbb{Z})$ described in \cite{Freitag} and elsewhere.

\ One can readily adapt arguments given in  \cite{Freitag} to see that the group $\Gamma_{p}$ acts totally discontinuously on $H^{+}(\mathbb{R}^{n})$. As $\Gamma_{p}[N]\subset\Gamma_{p}$ for each positive integer $N$ it follows that $\Gamma_{p}[N]$ also acts totally discontinuously on $H^{+}(\mathbb{R}^{n})$. Consequently $H^{+}(\mathbb{R}^{n})/\Gamma_{p}[N]$ is a conformally flat manifold. We will denote this manifold by ${\cal{M}}_{p}[N]$.
 
\subsection{Fundamental domains and their associated conformally flat manifolds}

It is known that all discrete arithmetic subgroups $\Gamma$ of Ahlfors-Vahlen's group $SAV(\mathbb{R}^{n-1})$ possess a fundamental domain in
$H^+(\mathbb{R}^n)$. See for example \cite{EGM88} where the general $n$-dimensional case is treated in condensed form.  A very detailed description of the fundamental domains in the particular three-dimensional case can be found in \cite{EGM98}. Just to recall: 
For a discrete subgroup, $\Gamma$, of $SAV(\mathbb{R}^{n-1})$ a fundamental domain ${\cal{F}}(\Gamma) \subset H^+(\mathbb{R}^n)$ is a relatively closed domain in $H^+(\mathbb{R}^n)$ with the two properties:
\par\medskip\par
1. $H^+(\mathbb{R}^n) = \cup_{M \in \Gamma} M\langle {\cal{F}}(\Gamma) \rangle$ \\
2. $int\;{\cal{F}}(\Gamma) \cap M \langle int\;{\cal{F}}(\Gamma)\rangle \neq 0, M \in \Gamma, \Longrightarrow  M = \pm I$.

\par\medskip\par

One can describe the geometry of a fundamental domain of a discrete group
$\Gamma$ acting on the upper half space in $\mathbb{C}$ in terms of the set of its isometric
circles
\cite{bulla,Ford}.

\par\medskip\par

The isometric circle for a M\"{o}bius transformation $\frac{az+b}{cz+d}$ where $\left(\begin{array}{cc} a & b \\ c & d \end{array}\right) \in SL(2,\mathbb{Z})$ with $c \neq 0$ would be the circle
$$
\{z \in H^+(\mathbb{C}):\|cz+d\|=1\}.
$$
In \cite{bulla} it is shown that the role played by isometric circles in introducing fundamental
domains for the arithmetic subgroups of the modular group $SL(2,\mathbb{Z})$, also referred to as Ford domains, can be carried over to upper half space in $\mathbb{R}^{n}$.
In this context the isometric circles are replaced by isometric spheres. So given a M\"{o}bius transformation $(ax+b)(cx+d)^{-1}$ induced by a matrix $M = \left(\begin{array}{cc} a & b \\ c & d \end{array}\right) \in SAV(\mathbb{R}^{n})$ with $c\ne 0$ we have that 
$$
\Big\|\frac{\partial [(ax+b)(cx+d)^{-1}]}{\partial x_{j}}\Big\|=\frac{1}{\|cx+d\|^{2}} \quad \mbox{for}\quad j=1,\ldots,n.
$$
Consequently for such a M\"{o}bius transformation its
isometric sphere is defined to be the sphere $\{x\in\mathbb{R}^{n}:\|cx+d\|=1\}$. This is a sphere in $\mathbb{R}^{n}$ centered at $-dc^{-1}$ and of radius $\|dc^{-1}\|$. Let us
denote this isometric sphere by $S_{M}$. Following \cite{bulla,zoll} in complete parallel to the complex case it may be shown that $M<S_{M}>=S_{M^{-1}}$. In analogy to the complex case we say that the
M\"{o}bius transformation $(ax+b)(cx+d)^{-1}$ induced by $M = \left(\begin{array}{cc} a & b \\ c & d \end{array}\right) \in SAV(\mathbb{R}^{n})$ with $c\neq 0$ is hyperbolic if the isometric  spheres $S_{M}$ and $S_{M^{-1}}$ are external. It is
elliptic if they intersect and $M<x>$ is parabolic if $S_{M}$ and $S_{M^{-1}}$ are tangent. It should be noted that as in the case of isometric circles the interior
$B(-dc^{-1},\|dc^{-1}\|)$, of the isometric sphere $S_{M}$ is mapped via $M$ to the exterior of the isometric sphere $S_{M^{-1}}$.

\par\medskip\par

Suppose $H$ is a discrete subgroup of $SAV(\mathbb{R}^{n-1})$ that acts totally discontinuously on $\mathbb{R}^{n}$, and ${\cal{F}}(H)$ is its associated fundamental domain. If one glues exactly the equivalence points of the boundary parts $\partial {\cal{F}}_p$ under the group action together, then one obtains a conformally flat manifold. This results from the factorization
$H^+(\mathbb{R}^n)/H$.
These manifolds belong to the general class of hyperbolic manifolds. 

\par\medskip\par

Let us give some simple examples.

\par\medskip\par 
 
For illustration, let us consider as a concrete example the group $\Gamma_{1}[4]$. This group is isomorphic to the classical arithmetic group $\Gamma[4]$. In \cite{bulla, fk} it is shown that one fundamental
domain for $H^+(\mathbb{C})/\Gamma[4]$ is the open set in $H^+(\mathbb{C})$ bounded between the lines $x=-1$ and $x=3$ and the $8$ isometric circles of radius $\frac{1}{4}$
and centered on the $x$ axis at  $\frac{-3}{4},\ldots,2\frac{3}{4}$ respectively. These are the isometric circles for the M\"{o}bius transformations
$$
M_{1}<z>\;=\;\frac{5z+4}{-4z-3},\ldots,M_{8}<z> \;=\; \frac{-3z+8}{4z-11},\;\; M_1,\ldots,M_8 \in\Gamma[4],
$$
respectively. To obtain the corresponding Riemann surface one identifies the two lines $\{z=x+iy:x=-1\}$ and
$\{z=x+iy:x=3\}$ and the isometric semicircle lying in upper half space associated with $M_{1}<z>$ with the isometric semicircle associated with $M_{8}<z>$. Further for  $j=2$,
$4$, $6$ the isometric semicircles lying in upper half plane and associated with the M\"{o}bius transformation $M_{j}<z>$  is identified with the isometric semicircle lying in
upper half space and associated with $M_{j+1}<z>$. If we consider the action of $\Gamma_{1}[4]$ on $H^+(\mathbb{R}^n)$, then the isometric semicircles associated with $M_{1},\ldots,M_{8}$ are replaced by the
isometric hemispheres $C_{1},\ldots,C_{8}$ lying in upper half space $H^+(\mathbb{R}^n)$ and now associated with the M\"{o}bius transformations
$$
M'_{1}<x>=(5x+4e_{1})(4e_{1}x-3)^{-1},\ldots,M'_{8}<x>=(-3x+8e_{1})(-4e_{1}x-11)^{-1}.
$$
Furthermore the lines $x=-1$ and $x=3$ are replaced by the half hyperplanes
$C_{9}=\{x\in H^+(\mathbb{R}^n):x_{1}=-1\}$ and $C_{10}=\{x\in H^+(\mathbb{R}^n):x_{1}=3\}$. A fundamental domain associated to the action of
$\Gamma_{1}[4]$ on $H^{+}(\mathbb{R}^n)$ is the domain in
$H^+(\mathbb{R}^n)$ with boundary $C_{1}\cup\ldots\cup C_{10}$. This domain is unbounded in the variables $x_{2},\ldots,x_{n-1}$ and $x_{n}$. Consequently the boundary of this fundamental domain has nontrivial intersection with unbounded open subsets of $\partial H^{+}(\mathbb{R}^{n})$. 

\ To obtain the conformally flat manifold $H^+(\mathbb{R}^n)/\Gamma_{1}[4]$ we now identify $C_{9}$ with
$C_{10}$, $C_{1}$ with $C_{8}$ and $C_{j}$ with $C_{j+1}$ for $j=2$, $4$, $6$.

\par\medskip\par

The group $\Gamma_2$ is isomorphic to the special linear group over the Gaussian integers $SL_2(\mathbb{Z}[e_1])$, where $\mathbb{Z}[e_{1}]$ is the lattice of Gaussian integers $\{p+qe_{1}:p,q\in\mathbb{Z}\}$. This group acts totally discontinuously on the modified upper half-space 
$$H^{+}(\mathbb{R} \oplus \mathbb{R}^{n-1})
= \{x_0 + x_1 e_1 + \cdots + x_{n-1}
e_{n-1}\;:
x_{n-1}
> 0\}.$$ Of course the group $SL_{2}(\mathbb{Z}[e_{1}])$ is essentially the Picard group $SL_{2}(\mathbb{Z}[i])$, the subgroup of $SL_{2}(\mathbb{C})$ with coefficents in the lattice of Gaussian integers $\{p+qi:p,q\in\mathbb{Z}\}$. Following for example \cite{Kri88}, one may determine that one fundamental domain for $SL_2(\mathbb{Z}[e_1])$ is  
$F_{2} := \{x=x_{0}+x_{1}e_{1}+\ldots+x_{n-1}e_{n-1} \in H^+(\mathbb{R}\oplus \mathbb{R}^{n-1})$ $\;|\; \|x\| \ge 1,\;0 < x_j < \frac{1}{2},\;j=0,1\}.$

\ It follows from our constructions that standard methods of constructing fundamental domains for $\Gamma_{1}[N]$ in one complex variable and for $\Gamma_{2}[N]$ in three real variables described in \cite{Ford,fk,Kri88} and elsewhere extend readily to give us examples of fundamental domains  for $\Gamma_{1}[N]$ and $\Gamma_{2}[N]$ in $n$-real variables.
 
\ Notice that for
a given discrete subgroup $\Gamma_{p}[N]$ there are infinitely many choices of fundamental domain. If $U$ is a particular fundamental domain for $\Gamma_{p}[N]$ then so is $M<U>$ for each $M\in\Gamma_{p}[N]$. The fundamental domain we will choose to work with will be denoted by ${\cal{F}}_{p}[N]$ and it will be a fundamental domain for $\Gamma_{p}[N]$ that lies in the set $U(p,N)=\{x\in H^{+}(\mathbb{R}^{n}):\frac{-N}{2}\le x_{i} \le \frac{N}{2}$ for $i=1,\ldots, p\}$, and it will be unbounded in the variables $x_{p+1},\ldots, x_{n}$.

\section{Monogenic Eisenstein series on $\Gamma_p[N]$ and the Hecke trick}

$H^{+}(\mathbb{R}^n)$ is the universal covering space of the class of manifolds ${\cal{M}}_{p}[N]$ we are considering here.  As a consequence, there
exists a well-defined projection map $p: H^{+}(\mathbb{R}^n) \rightarrow {\cal{M}}_{p}[N]: x \mapsto x\;(mod\; \Gamma_p[N])$. Let us write $x'$ for $p(x)$ where $x \in H^{+}(\mathbb{R}^{n})$. For each open set $U\subseteq H^{+}(\mathbb{R}^n)$ we write $U'$ for $p(U)$ which in turn is an open subset of ${\cal{M}}_{p}[N]$. Further for each set $S\subset {\cal{F}}_{p}[N]$ we will write $S'$ for $p(S)$.

\par\medskip\par

In this section we give some elementary examples of automorphic forms on the groups $\Gamma_p[N]$ (with $N \ge 3$). More precisely, the examples we give here
are invariant under the action of $\Gamma_p[N]$ up to the particular automorphic weight factor $J(M,x)$.  These then project down to form non-trivial examples of monogenic
sections on the associated families of conformally flat
 manifolds.
These sections take values in a fixed spinor bundle $F_1$. 

\par\medskip\par

More precisely this bundle is
constructed over ${\cal{M}}_{p}[N]$ by identifying each pair $(x,X)$ with $(M<x>,J(M,x)X)$ for every
$M \in \Gamma_p[N]$ where
$x \in H^+(\mathbb{R}^n)$, and $ X \in Cl_n$.

\par\medskip\par

It should be noted that as $(ax+b)(cx+d)^{-1}=(-ax-b)(-cx-d)^{-1}$ then usually there is an ambiguity of sign in the previous identification of $X$ with $J(M,x)X$. However as $-I$ does not belong to the group $\Gamma_{p}[N]$ for $N>2$ there is no ambiguity of sign in this particular case and so the bundle $F_{1}$ is globally well defined.

\begin{definition}
Suppose that $U'$ is a domain in ${\cal{M}}_{p}[N]$
and $f:U'\rightarrow F_{1}$ is a section that locally is annihilated on the left by the Dirac operator, $D$, then $f$ is called a (left) monogenic section.
\end{definition}

In the previous definition we locally used the Dirac operator, $D$. In fact we have introduced a Dirac operator
$D_{{\cal{M}}_{p}[N]}$
that acts globally on sections in $F_{1}$. This Dirac operator is in fact the Atiyah-Singer operator, or Atiyah-Singer-Dirac operator, cf. \cite{LM}. 

\par\medskip\par

\ In \cite{KraHabil} it is noted that when $N=1$ or $N=2$ the matrix $-I$ belongs to $\Gamma_{N}$. Consequently any function that satisfies $f(x) = J(M,x) f(M\langle x \rangle)$ for all $M \in \Gamma_p[N]$ and each $x\in H^{+}(\mathbb{R}^{n})$ must satisfy $f(x)=-f(x)$ and so vanishes identically. For this reason we shall unless otherwise specified work in this paper with the cases $N\geq 3$.

\par\medskip\par

For $N \ge 3$, as $-I$ is not in $\Gamma_p[N]$ it is possible to construct non-trivial monogenic functions that satisfy $f(x)=J(M,x)f(M,\langle x\rangle)$ for all $M\in\Gamma_{p}[N]$. To introduce these types of Eisenstein series we shall need the following convergence lemma which is proved in \cite{KraHabil}.

\begin{lemma}
(Convergence lemma) For all $\alpha>p+1$ and all $N\in\mathbb{N}$ the series
$$
\sum\limits_{M:\Gamma_{p}[N]\backslash{\cal{T}}_{p}[N]}\frac{1}{\|cx+d\|^{\alpha}}
$$
converges uniformly on each compact subset of $H^{+}(\mathbb{R})$. Here $M:\Gamma_{p}[N]\backslash{\cal{T}}_{p}[N]$ denotes a sum over representatives of the left coset space $\Gamma_{p}[N]\backslash{\cal{T}}_{p}[N]$.
\end{lemma}

\ Now following \cite{KraHabil} we introduce the following generalized Eisenstein series.

\begin{definition}
For $p < n-2$ the monogenic Eisenstein series attached to $\Gamma_p[N]$ are then defined by
\begin{equation}
\label{eisenstein}
E_{p,N}(x) = \sum\limits_{M:\Gamma_p[N]\backslash{\cal{T}}_{p}[N]} J(M,x).
\end{equation}
\end{definition}
Here the notation $M:\Gamma_p[N]\backslash{\cal{T}}_{p}[N]$ means that the matrices $M$ run through a system of representatives of right cosets in $\Gamma_p[N]$ modulo the translation group ${\cal{T}}_{p}[N]$.  
\ The convergence lemma insures the convergence of the Series (\ref{eisenstein}). 

\ Here we shall assume that each representative $M$ from the left coset space $\Gamma_{p}[N]\backslash{\cal{T}}_{p}[N]$ is chosen so that $M({\cal{F}}_{p}[N]$ lies in the set $U(p,N)$.

\ Since the identity matrix is a representative from the coset space $\Gamma_{p}[N]\backslash{\cal{T}}_{p}[N]$ and $J(I,x)=1$ and as  $\lim_{x\rightarrow\infty}J(M,x)=0$ for each $M\in\Gamma_{p}[N]\backslash{\cal{T}}_{p}[N]$ with $M\ne I$ then $\lim_{x_n \rightarrow +\infty} E_{p,N}(e_n x_n) = 1$. Consequently the Series (\ref{eisenstein}) does not vanish identically when $N\geq 3$. 

\ One can directly verify by a rearrangement argument, \cite{KraHabil}, that the Series (\ref{eisenstein}) satisfies the transformation rule
$$
E_{p,N}(x) = J(M,x) E_{p,N}(M<x>) \quad\quad \forall M \in \Gamma_p[N].
$$
Hence the generalized Eisenstein series (\ref{eisenstein}) project down to non-trivial 
monogenic sections with values in the spinor bundle $F_1$ over the manifold ${\cal{M}}_{p}[N]$.

\par\medskip\par

{\bf Remark.} If we extend the sum in Expression (\ref{eisenstein}) to the whole group $\Gamma_p[N]$, the series
would diverge. This is due to the fact that the summation of the expressions $J(M,x)$ over the translation
group ${\cal{T}}_p[N]$ diverges.
\par\medskip\par
 By adapting the so called Hecke trick we can get convergent Eisenstein series in the remaining two cases $p=n-2$ and $p=n-1$. This is a classical
method to introduce Eisenstein and Poincar\'e series of lower weight. See for instance \cite{Freitag}.

In analogy to the one complex variable and several complex variable cases let us introduce the following adapted Eisenstein series
\begin{equation}
E_{p,N}(x)(x,s) = \sum\limits_{M: {\mathcal{T}}_p[N] \backslash \Gamma_p[N]} \Bigg(\frac{x_n}{\|cx+d\|^2}\Bigg)^s J(M,x),
\end{equation}
where $s$ is a complex auxiliary parameter.
In view of the convergence lemma these
series are a priori normally convergent and complex-analytic in $s$ 
whenever the real part of $s$ is greater than  $- n+2+p$. 

\ Whenever the real part of $s$ is greater than  $-n+2+p$, the expressions $\Bigg(\frac{x_n}{\|cx+d\|^2}\Bigg)^s J(M,x)$ are not monogenic in the vector variable $x$. However, one can show:
\begin{proposition}\label{prop2}
The series $E_{p,N}(x,s)$ possesses a continuous extension to the complete complex $s$-semiplane $\{s=u+iv\in\mathbb{C}:v\geq 0\}$.
\end{proposition}
The limit 
\begin{equation}
E_{p,N}(x) := \lim\limits_{s \rightarrow 0^+}\sum\limits_{M:{\mathcal{T}}_p[N] \backslash \Gamma_p[N]}  \Bigg(\frac{x_n}{\|cx+d\|^2}\Bigg)^s J(M,x)
\end{equation}
will then provides us with a left monogenic Eisenstein series for the larger groups $\Gamma_{n-2}[N]$ and $\Gamma_{n-1}[N]$ which then of course has the desired invariance behavior
under the respective actions of  $\Gamma_{n-2}[N]$ and $\Gamma_{n-1}[N]$.

\par\medskip\par

The proof of the continuous extension of the series $E_{p,N}(x,s)$ towards $s \rightarrow 0^+$ can be done in the same way as  the classical proof of the holomorphic Hilbert
Eisenstein series in several complex variables as presented for instance in \cite{Freitag} pp. 165--172. Hence, we leave the detailed proof as an exercise for the reader and restrict
ourselves to only present here the features of the proof which are different from the classical proof.

\par\medskip\par

Without loss of generality we focus on the case $p=n-1$, since the series  $E_{n-2,N}(x,s)$  are subseries of  $E_{n-1,N}(x,s)$.
As in the classical case, one expands the $\Gamma_{n-1}[N]$-periodic function $E_{n-1,N}(x,s)$ into a Fourier series. This can be done easily by first expanding the $n-1$-fold
periodic subseries
\begin{equation}
\varepsilon_{n-1,N}(x,s) = \sum\limits_{\underline{g} \in N \mathbb{Z}^{n-1}} G(x+\underline{g}) \Bigg(\frac{x_n}{\|x+ \underline{g}\|^2}\Bigg)^{s}
\end{equation}
into a Fourier series of the form
$$
\sum\limits_{\underline{g} \in \frac{1}{N}\mathbb{Z}^{n-1}} \alpha(x_n,\underline{g};s) e^{2 \pi i<\underline{g},\underline{x}>}
$$
and showing that this series has a continuous extension as $s\rightarrow 0^+$.
The Fourier coefficients are basically given (up to a constant) by the following integral (putting $\underline{x}:=x_1 e_1 + \cdots x_{n-1} e_{n-1}$):
\begin{equation}
\label{fourriercoeff}
\alpha(x_n,\underline{g};s) = \int\limits_{\mathbb{R}^{n-1}} D_y\Bigg[ \frac{x_n^s}{\|x+y\|^{n-2+2s}} e^{- 2 \pi i<x+y,\underline{g}>}\Bigg] d\underline{x}.
\end{equation}
We consider
\begin{eqnarray}
&=& \int\limits_{\mathbb{R}^{n-1}} \frac{x_n^{s}}{(\|\underline{x}\|^2+x_n^2)^{\frac{n-2}{2}+s}} e^{-2 \pi i <\underline{x},\underline{g}>} d\underline{x}\nonumber \\
&=& \int\limits_{S_{n-3}}\int\limits_{\theta=0}^{2\pi} \int\limits_{r=0}^{\infty} \frac{x_n^s r^{n-2}}{(r^2+x_n^2)^{\frac{n-2}{2}+s}} e^{-2 \pi i r \|\underline{g}\| \cos(\theta)} \sin^{n-3}(\theta) dr d\theta dS_{n-3}\nonumber \\
&=& \omega_{n-3}\int\limits_{r=0}^{\infty} \frac{x_n^s r^{n-2}}{(r^2+x_n^2)^{\frac{n-2}{2}+s}}r^{-\frac{n-3}{2}} J_{\frac{n-3}{2}}(2 \pi r \|\underline{g}\|)dr \label{superint}\\
&=& \frac{2^{\frac{3-2s}{2}} (2 \pi |\underline{g}|)^{\frac{2s-1}{2}} \pi^{\frac{n-1}{2}} x_n^{\frac{1}{2}}}{\Gamma(\frac{n+2s}{2}-1)} K_{\frac{1-2s}{2}}(2 \pi \|\underline{g}\|x_n) 
\end{eqnarray}
where we applied the substitution $r :=\|\underline{x}\|$ and where $\theta$ is the angle between $\underline{x}$ and $\underline{g}$. $J_{n-3/2}$ and $K_{1-2s/2}$ are the standard Bessel functions of first and second type.

 The function $x_n^{\frac{1}{2}} K_{\frac{1-2s}{2}}$ is entire-complex analytic in $s$.  In particular, it is continuous at $s=0$. This  establishes the existence of the limit $s\rightarrow 0^+$. Therefore, we may in particular interchange the application of the Dirac operator in (9) with the integration process.

\par\medskip\par

The rest of the analyticity proof can now be adapted directly from the classical proof presented in \cite{Freitag} pp. 165. One re-expresses the complete series $E_{n-1,N}(x)$ as
a series over the Fourier series of $\varepsilon_{n-1,N}(x,s)$. 

\ As the term $\lim_{s\rightarrow 0}\sum\limits_{M_{1},\ldots, M_{r}}(\frac{x_{n}}{\|cx+d\|})^{s}J(M,x)$ is left monogenic, where summation is taken over any finite subset $\{M_{1},\ldots,M_{r}\}$ of $\Gamma_{k}[N]\backslash{\cal{T}}_{k}[N]$, it follows that the functions $E_{n-2,N}(x)$ and $E_{n-1,N}(x)$ are left monogenic.

\par\medskip\par

The projections of the functions $E_{n-2,N}(x)$ and $E_{n-1,N}(x)$ to sections over ${\cal{M}}_{p}[N]$ thus provide us also with non-trivial examples of
monogenic sections on the manifolds
${\cal{M}}_{p}[N]$ with $p=n-2$ and $p=n-1$.

\ We shall need the following, cf. \cite{KraHabil}:

\begin{definition}
Suppose for $k$ odd that $N\geq 3$ and for $k$ even $N\in\mathbb{N}$. Suppose also that $p<$ then the $k$-genic Eisenstein series attached to $\Gamma_{p}[N]$ are defined to be the series
\begin{equation}
\label{Eisenk}
E_{p,N,k}(x)=\sum\limits_{M:\Gamma_{p}[N]\backslash{\cal{T}}_{p}[N]}J_{k}(M,x).
\end{equation}
\end{definition}

\ Convergence of the Series (\ref{Eisenk}) follows from Lemma 2.

\section{Cauchy kernels for monogenic and $k$-genic sections}

\ Here we use the Eisenstein series introduced in Definition 3 to introduce an explicit formula for the Cauchy kernel or fundamental solution to the Dirac operator on the manifolds ${\cal{M}}_{p}[N]$
with $p < n-2$. We then use the Hecke trick introduced in the previous section to introduce the fundamental solutions for the cases $p=n-2$ and $p=n-1$. We also introduce fundamental solutions to the analogues of the operators $D^{k}$ for $k<n$ over ${\cal{M}}_{p}[N]$ using the Eisenstein series introduced in Definition 4. Further Calderon-Zygmund operators in this context are introduced.

\par\medskip\par

We now proceed to prove the first main result of this section:

\begin{proposition}
For $p < n-2$ and for each $x\in{\cal{F}}_{p}[N]$ the series
\begin{equation}
\sum\limits_{T\in{\cal{T}}_{p}[N]}\sum\limits_{M:\Gamma_{p}[N]\backslash{\cal{T}}_{p}[N]}J(TM,x)G(y-TM<x>).
\end{equation}
converges uniformly on any compact subset $K$ of ${\cal{F}}_{p}[N]\backslash\{x\}$.
\end{proposition}
{\bf Remark.} Notice that the second series $\sum\limits_{M:\Gamma_{p}[N]\backslash{\cal{T}}_{p}[N]}J(TM,x)G(y-TM<x>)$ depends on the particular choice of the system of representatives of right cosets in $\Gamma_p[N]$ modulo ${\cal{T}}_p[N]$ since the function $G(y-TM<x>)$ is not invariant under  ${\cal{T}}_{p}[N]$. This notation is understood in the sense that one has to specify first a particular system of representatives. The whole double sum however then gets again independent of that particular choice since the first sum extends ober the whole translation group ${\cal{T}}_p[N]$.  
\par\medskip\par
{\bf{Proof}}: It should first be noted that $J(TM,x)=J(M,x)$ for each $T\in{\cal{T}}_{p}[N]$ and each $M\in\Gamma_{p}[N]$. So from Lemma 1 for each fixed $T\in{\cal{T}}_{p}[N]$ the series $\sum\limits_{M\in\Gamma_{p}[N]\backslash{\cal{T}}_{p}[N]}J(TM,x)$ converges uniformly on $K$ for $p<n-2$.  Further the series (13) can be rewritten as
$$
\sum\limits_{T\in{\cal{T}}_{p}[N]}\sum\limits_{M:\Gamma_{p}[N]\backslash{\cal{T}}_{p}[N]}J(M,x)G(y-TM<x>).
$$
\ In \cite{KraHabil} it is shown that the series $\sum\limits_{T\in{\cal{T}}_{p}[N]}G(y-T<x>)$ converges uniformly for $p<n-1$. As each representative $M$ from $M\in\Gamma_{p}[N]$ is chosen so that $M({\cal{F}}_{p}[N])\subset U(p,N)$ then for each $T_{\underline{m}}\in{\cal{T}}_{p}[N]\backslash\{I\}$ and each $y\in{\cal{F}}_{p}[N]$ we have that $\|G(y-T_{\underline{m}}M<x>)\|\leq \|G(w-\underline{m})\|$ where $T_{\underline{m}}<x>=x+\underline{m}$ and $\underline{m}\in\mathbb{Z}e_{1}\ldots\mathbb{Z}e_{p}$. Consequently it follows from Lemma 1 that the subseries 
$$
\sum\limits_{T\in{\cal{T}}_{p}[N]\backslash\{I\}}\sum\limits_{M:\Gamma_{p}[N]\backslash{\cal{T}}_{p}[N]}J(TM,x)G(y-TM<x>)
$$ 
is uniformly convergent on $K$ for $p<n-2$. 

\ Further for $\{B(x_{1},r_{1},\ldots , B(x_{q},r_{q})\}$ a finite covering of $K$ by open balls whose closure lie in ${\cal{F}}_{p}[N]\backslash\{x\}$ let $r(K)$ denote the minimum radius of these balls. Then for each $M\in\Gamma_{p}[N]$ we have that $\|G(y-M<x>)\|\leq \frac{1}{\omega_{n}}r(K)^{1-n}$. It now follows that the series $\Sigma_{M:\Gamma_{p}[N]\backslash{\cal{T}}_{p}[N]}\|J(M,x)G(y-M<x>)\|$ is dominated by the series $\frac{1}{\omega_{n}}r(K)^{1-n}\Sigma_{M:\Gamma_{p}[N]\backslash{\cal{T}}_{p}[N]}\|J(M,x)\|$. Consequently it follows from Lemma 2 that the subseries $\Sigma_{M:\Gamma_{p}[N]\backslash{\cal{T}}_{p}[N]}J(M,x)G(y-M<x>)$ is convergent on $K$ for $p<n-2$. 

\ It now follows that the series (13) is uniformly convergent on $K$ for $p<n-2$. $\Box$

\ In fact the above proof can be readily adapted to show that the above series is uniformly convergent on any compact subset of $H^{+}(\mathbb{R}^{n})\backslash\cup_{M\in\Gamma_{p}[N]}\{M<x>\}$. Moreover it may be readily seen that this series defines a right monogenic function, $C_{p,N}(x,y)$,  in the variable $y$ on the domain $H^{+}(\mathbb{R}^{n})\backslash\cup_{M\in\Gamma_{p}[N]}\{M<x>\}$. Further as the function $J(M,x)G(w-M<x>)$ is left monogenic in $x$ the function $C_{p,N}(w,x)$ is left monogenic in $x$.

\ It should be noted that it follows from the definition of a fundamental domain that the kernel $C_{p,N}(x,y)$ has precisely one singularity in the domain ${\cal{F}}_{p}[N]$.This singularity is of order $n-1$.

\begin{proposition}
The kernel $C_{p,N}(x,y)$ satisfies the asymmetry relation $C_{p,N}(x,y)=-\tilde{C}_{p,N}(y,x)$.
\end{proposition} 
{\bf{Proof}}:
Notice that
$$
G(y-M<x>) =G(y-(ax+b)(cx+d)^{-1})=G(y-(\widetilde{cx+d})^{-1}(\widetilde{ax+b}))
$$
$$
=G(y-(x\tilde{c}+\tilde{d})^{-1}(x\tilde{a}+\tilde{b}))=(x\tilde{c}+\tilde{d})^{-1}\|x\tilde{c}+\tilde{d}\|^{n}
\frac{
(x \tilde{c}+\tilde{d}) y-x\tilde{a}-\tilde{b}
}
{
\|(x\tilde{c}+\tilde{d})y-x\tilde{a}-\tilde{b}\|^{n}
}
$$
$$
=(\widetilde{cx+d})^{-1}\|cx+d\|^{n}\frac{
(x\tilde{c}+\tilde{d})y-x\tilde{a}-\tilde{b}
}
{
\|(x\tilde{c}+\tilde{d})y-x\tilde{a}-\tilde{b}\|^{n}
}.
$$
So $J(M,x)G(y-M<x>)=\frac{(x\tilde{c}+\tilde{d})y-x\tilde{a}-\tilde{b}}{\|(x\tilde{c}+\tilde{d})y-x\tilde{a}-\tilde{b}\|^{n}}$.

\ Now
$$
(x\tilde{c}+\tilde{d})y-x\tilde{a}-\tilde{b}=x\tilde{c}y+\tilde{d}y-x\tilde{a}-\tilde{b}
$$
$$
=x(\tilde{c}y-\tilde{a})+\tilde{d}-\tilde{b}=-x(-\tilde{c}+\tilde{a})+\tilde{d}y-\tilde{b}
$$
$$
=(-x+(\tilde{d}y-\tilde{b})(-\tilde{c}y+\tilde{a})^{-1})(-\tilde{c}y+\tilde{a})=(-x+M^{\star}<y>)(\widetilde{-yc+a})
$$
where $M^{\star} = \left( \begin{array}{cc} \tilde{d} & -\tilde{b} \\ -\tilde{c} & \tilde{a}\end{array}\right)$ $\in\Gamma_{p}[N]$.  Consequently
$$
J(M,x)G(y-M<x>)=G(M^{\star}<y>-x)\tilde{J}(M^{\star},y).
$$
As $\star:\Gamma_{p}[N]\rightarrow\Gamma_{p}[N]:M \rightarrow M^{\star}$ is an isomorphism it follows that $C_{p,N}(x,y)=-\tilde{C}_{p,N}(y,x)$. $\Box$

\ A further property of the kernel $C_{p,N}(x,y)$ is given as follows.

\begin{proposition}
The kernel, $C_{p,N}(x,y)$, satisfies the relation $J(L,x)C_{p,N}(L<x>,y)=C_{p,N}(x,y)$ for each $L\in\Gamma_{p}[N]$.
\end{proposition} 
{\bf{Proof}}:
Now
$$
J(M,L<x>) G(y-(ML)<x>)=J(L,x)^{-1}  J(ML,x) G(y-(ML)<x>)
$$
$$= J(L,x)^{-1}J(A,x)G(y-A<x>)
$$
where $A=ML$. As $\Theta_{L}:\Gamma_{p}[N]\rightarrow\Gamma_{p}[N]:M\rightarrow ML$ is a bijection we now have that
$C_{p,N}(L<x>,y)= J(L,x)^{-1} C_{p,N}(x,y)$. $\Box$

The projection $C'_{p,N}(x',y')$ of the kernel $C_{p,N}(x,y)$ is thus an $F_{1}$ valued left monogenic section in $y'$ that lives over the manifold ${\cal{M}}_{p}[N]$.
\par\medskip\par
We may therefore draw the conclusion that $C_{p,N}$ is a monogenic $\Gamma_p[N]$-periodic function. Its projection is a section on ${\cal{M}}_{p}[N]$ with only one point singularity. This singularity is of the order of the Cauchy kernel. 

\begin{definition}
A hypersurface $S$ lying in $\mathbb{R}^{n}$ is called a strongly Lipschitz surface if locally it is the graph of a real valued Lipschitz continuous function and the Lipschitz constants for these local Lipschitz graphs are bounded.
\end{definition}

\begin{theorem}
\label{cauchythm}
(Cauchy's integral formula for ${\cal{M}}_{p}[N]$)

Suppose $p<n-2$. Let $U \subset F_{p}[N]$ be a domain and $V$ be a bounded subdomain with a strongly Lipschitz boundary, $S$, lying in $U$. Then for each $y' \in V'$ and every left monogenic section
$f':U'\rightarrow F_1$ we have
\begin{equation}\label{cauchyint}
f'(y') =  \int\limits_{\partial S'} \tilde{C}'_{p,N}(x',y') d\sigma'(x') f'(x'),
\end{equation}
\end{theorem}
{\bf{Proof}}: We can rewrite the kernel $C_{p,N}(x,y)$ as $G(x-y)+(C_{p,N}(x,y)-G(x-y))$. The term $C_{p,N}(x,y)-G(x-y)$ has no singularity on the fundamental domain ${\cal{F}}_{p}[N]$ while $DG(x-y)=\delta_{x=y}$, the Dirac delta function. The result follows on projecting to the manifold ${\cal{M}}_{p}[N]$. $\Box$

\ This property has the following interpretation within the theory of automorphic forms:

Recall that all the sections that are well-defined  and monogenic on the whole manifold ${\cal{M}}_{p}[N]$ lift to the class of automorphic forms on $\Gamma_p[N]$  that have the
property of being monogenic on the whole upper half-space. In particular, it locally reproduces the projections to ${\cal{M}}_{p}[N]$ of the special monogenic Eisenstein series $E_{p,N}(x)$. Furthermore, it reproduces all the projections to ${\cal{M}}_{p}[N]$ of the general Poincar\'e series given in \cite{KraHabil} of the type
$$
Q(x,\tilde{f}) = \sum\limits_{M:\Gamma_p[N]\backslash{\cal{T}}_{p}[N]} J(M,x) \tilde{f}(M<x>).
$$
where $\tilde{f}$ is an arbitrary bounded left monogenic function on $H^{+}(\mathbb{R}^n)$ that is invariant under the translation group ${\cal{T}}_p[N]$.

\ For $k < n$, Lemma 1 and the proof of Proposition 3 can easily be adapted to also obtain explicit formulas for  kernels for $k$-genic sections on the manifold ${\cal{M}}_{p}[N]$. 

\ More precisely we have
\begin{proposition}
For $k < n$ and $p < n-1-k$ the series
\begin{equation}\label{ckg}
\sum\limits_{T\in{\cal{T}}_{p}[N]}\sum\limits_{M :\Gamma_p[N]\backslash{\cal{T}}_{p}[N]} J_k(TM,x)G_k(y-TM<x>).
\end{equation}
is uniformly convergent on any compact subset of ${\cal{F}}_{p}[N]\backslash\{x\}$.
\end{proposition}

\ We shall denote the kernel given by expression (15) by $C_{k,p,N}(x,y)$. Note that when $k=1$ we get back the kernel $C_{p,N}(x,y)$.

\ In the cases where $k$ is even the kernel $C_{k,p,N}(x,y)$ do not vanish in the particular cases where $N=1$ or $N=2$. This is because each term $J_{k}(M,x)G_{k}(M<x>-y)$ arising in Series (15) equal $\frac{1}{\|cx+d\|^{n-k}\|y-M<x>\|^{n-k}}$.  These are all positive terms so Series (15) cannot vanish.

\ By very similar arguments to those used to establish Proposition 4 it may be determined that when $k$ is odd and $N\geq 3$ then $C_{k,pN}(x,y)=-\tilde{C}_{k,p,N}(y,x)$ and when $k$ is even and $N\in\mathbb{N}$ then $C_{k,p,N}(x,y)=C_{k,p,N}(y,x)$. Further by a minor adaptation of arguments given to prove Proposition 5 it may be seen that $J_{k}(L,x)C_{k,p,N}(L<x>,y)=C_{k,p,N}(x,y)$ for each $L\in\Gamma_{p}[N]$.

\ We may now consider new bundles, $F_{k}$ constructed over ${\cal{M}}_{p}[N]$ by making the identification $(x,X)$ with $(M<x>,J_{k}(M,x)X)$ for each $M\in\Gamma_{p}[N]$ with $x\in H^{+}(\mathbb{R}^{n})$ and $X\in Cl_{n}$. Note that in the cases where $k$ is even the conformal weight factor $J_{k}(M,x)$ is real valued. So in these cases the bundles are not spinor bundles. 

We can now establish the following:
\begin{theorem}\label{gf} Suppose that $\psi':{\cal{M}}_{p}[N]\rightarrow F_{k}$ is a $C^{k}$ section with compact support. Then for each $y'\in{\cal{M}}_{p}[N]$
\begin{equation}
\psi'(y')=D_{{\cal{M}}_{p}[N]}^{k}\int_{{\cal{M}}_{p}[N]}\tilde{C}'_{k,p,N}(y',x')\psi'(x')dm'(x')
\end{equation}
and
\begin{equation}
\psi'(y')=\int_{{\cal{M}}_{p}[N]}\tilde{C}'_{k,p,N}(x',y')D_{{\cal{M}}_{p}[N]}^{k}\psi'(x')dm'(x')
\end{equation}
where $C'_{k,p,N}(x',y')$ is the projection to $F_{k}$ of the kernel $C_{k,p,N}(x,y)$ and $m'$ is the projection to ${\cal{M}}_{p}[N]$ of Lebesgue measure on ${\cal{F}}_{p}[n]$. Further $p<n-2$ and $N\geq 1$ when $k$ is even and $N\geq 3$ when $k$ is odd.
\end{theorem}
{\bf{Proof}}: The section $\psi'$ lifts to a function $\psi$ defined on the fundamental domain ${\cal{F}}_{p}[N]$. This function is $C^{k}$ and has compact support. Now consider $D^{k}\int_{{\cal{F}}_{p}[N]}\tilde{C}_{p,N}(x,y)\psi(x)dx^{n}$. This expression is equal to 
\begin{equation}
D^{k}\int_{{\cal{F}}_{p}[N]}G_{k}(x-y)\psi(x)dx^{n}+D^{k}I(y)
\end{equation}
where $I(y)=\int_{{\cal{F}}_{p}[N]}(\tilde{C}_{p,N}(x,y)\psi(x)-G_{k}(x-y)\psi(x))dx^{n}$. The term $I(y)$ is left $k$-genic on ${\cal{F}}_{p}[N]$. So expression (18) reduces to $D^{k}\int_{{\cal{F}}_{p}[N]}G_{k}(x-y)\psi(x)dx^{n}$ and this term is equal to $\psi(y)$.

\ If instead of lifting $\psi'$ to the fundamental domain ${\cal{F}}_{p}[N]$ we had lifted it to the fundamental domain $M<{\cal{F}}_{p}[N]>$ for some $M\in\Gamma_{p}[N]$ then $\psi'$ would lift to the function $\psi(u)$ where $u\in M<{\cal{F}}_{p}[N]>$. Again this function is $C^{k}$ with compact support in $M<{\cal{F}}_{p}[N]>$. In this case we would get by changing variables back to $x$ and $y$
$$
J_{-k}(M,y)\psi(M<y>)=D^{k}\int_{{\cal{F}}_{p}[N]}\tilde{C}_{p,N}(x,y)J_{-k}(M,x)\psi(x)dx^{n},
$$
where $J_{-k}(M,x)=\frac{\widetilde{cx+d}}{\|cx+d\|^{n+k+1}}$ when $k$ is odd and $J_{-k}(M,x)=\frac{1}{\|cx+d\|^{n+k}}$ when $k$ is even.
 This establishes equation (16). Similar arguments gives us equation (17). $\Box$

\ It remains to set up a monogenic Cauchy kernel for the cases $p=n-2, p=n-1$. On placing $p=n-1$ or $p=n-1$ in the series for $C_{p,N}(x,y)$ this series will now diverge. However, we can bypass this issue by again adapting the Hecke trick. For each $M\in\Gamma_{p}[N]$ let us write the term $\frac{x_{n}}{\|cx+d\|^{2}}$ as $H(M,x)$. It should be noted that for each $T\in{\cal{T}}_{p}[N]$ and each $M\in\Gamma_{p}[N]$ we have that $H(TM,x)=H(M,x)$. Now let us introduce the series
$$
\sum\limits_{T\in{\cal{T}}_{p}[N]}\sum\limits_{M:\Gamma_p[N]\backslash{\cal{T}}_{p}[N]} H(TM,x)^{s} J(TM,x) G(y-TM<x>),
$$
where $s$ is the complex auxiliary parameter introduced in Proposition ~\ref{prop2}. As $H(TM,x)=H(M,x)$ and $J(TM,x)=J(M,x)$ one can apply the Hecke trick described in Section 3  to adapt the proof of Proposition 3 to establish the following:

\begin{proposition}
The series 
\begin{equation}
\sum\limits_{T\in{\cal{T}}_{n-2}[N]}\sum\limits_{M:\Gamma_{n-2}\backslash{\cal{T}}_{n-2}[N]}H(TM,x,)^{s}J(TM,x)G(y-TM<x>)
\end{equation}
is absolutely convergent on any compact subset of ${\cal{F}}_{n-2}[N]\backslash\{y\}$ and for any $s\in\mathbb{C}$ whose real part is greater than zero.
\end{proposition}

\ Given that $H(TM,x)=H(M,x)$ and $J(TM,x)=J(M,x)$ the series (19) can be rewritten as 
\begin{equation}
\sum\limits_{T\in{\cal{T}}_{n-2}[N]}\sum\limits_{M:\Gamma_{n-2}\backslash{\cal{T}}_{n-2}[N]}H(M,x)^{s}J(M,x)G(y-TM<x>).
\end{equation}

\ Further one can combine the arguments used to establish Proposition 7 with arguments presented in \cite{KraHabil} to obtain

\begin{proposition}
The series
$$
\sum\limits_{M:\Gamma_{n-1}[N]\backslash{\cal{T}}_{p}[N]}H(M,x)^{s}J(M,x)G(y-M<x>)+\sum\limits_{T_{\underline{m}},T_{-\underline{m}}\in{\cal{T}}_{n-1}[N]}
$$
\begin{equation}
\sum\limits_{M:\Gamma_{n-1}[N]\backslash{\cal{T}}_{n-1}[N]}H(M,x)^{s}J(M,x)(G(y-T_{\underline{m}}M<x>)+G(y-T_{-\underline{m}}M<x>))
\end{equation}
is absolutely convergent on any compact subset of ${\cal{F}}_{n-1}[N]\backslash\{y\}$ and for any $s\in\mathbb{C}$ whose real part is positive.
\end{proposition}

\ Let us denote the series (20) and (21) by $C_{n-2,N,s}(x,y)$ and $C_{n-1,N,s}(x,y)$ respectively. From Proposition ~\ref{prop2} it now follows that $\lim_{s\rightarrow 0}C_{n-2,N,s}(x,y)$ and $\lim_{s\rightarrow 0}C_{n-1,N,s}(x,y)$ exist and are functions $C_{n-2,N}(x,y)$ and $C_{n-1,N}(x,y)$ respectively defined on $H^{+}(\mathbb{R}^{n})\times H^{+}(\mathbb{R}^{n})\backslash diagonal H^{+}(\mathbb{R}^{n})$ where $diagonal H^{+}(\mathbb{R}^{n})=\{(x,x):x\in H^{+}(\mathbb{R}^{n})$. 

\ Note that for each $L\in\Gamma_{p}[N]$ 
$$
H(M,L<x>)=H(ML,x).
$$
Consequently from the same arguments used to prove Proposition 5 we have that $J(L,x)C_{n-2,N,s}(L<x>,y)=C_{n-2,N,s}(x,y)$ and $J(L,x)C_{n-1,N,s}(L<x>,y)=C_{n-1,N,s}(x,y)$ for each $L\in\Gamma_{k}[N]$. It follows that 
$$
\lim_{s\rightarrow 0}C_{p,N,s}(L<x>,y)=J(L,x)^{-1}C_{p,N}(x,y)
$$ 
for $p=n-2$ and for $p=n-1$. So $J(L,x)C_{p,N}(L<x>,y)=C_{p,N}(x,y)$ also for $p=n-2$ and for $p=n-1$.

\ By the same argument used to prove Proposition 4 we may see that $C_{p,N,s}(y,x)=-\tilde{C}_{p,N,s}(x,y)$ for $p=n-2$ and $p=n-1$. Consequently $C_{p,N}(y,x)=-\tilde{C}_{p,N}(x,y)$ for $p=n-2$ and $p=n-1$.

\par\medskip\par
Notice that for $p=n-2$ and $p=n-1$ the functions $C_{p,N}(x,y)$ have exactly one point
singularity in each fundamental domain. This becomes clear when rewriting $C_{p,N}(x,y)$ in
the equivalent form $G(x-y)+(C_{p,N}(x,y)-G(x-y))$. The term $C_{p,N}(x,y)-G(x-y)$ has no singularities on the fundamental domain ${\cal{F}}_{p}[N]$. 

\ Bearing these comments in mind we can adapt arguments given in \cite{MM} and elsewhere to obtain:

\begin{theorem}
Suppose that $U$ is a bounded domain in ${\cal{F}}_{p}[N]$ with strongly Lipschitz boundary $S$ lying in ${\cal{F}}_{p}[N]$. Suppose further that $\psi:S\rightarrow Cl_{n}$ belongs to $L^{q}(S)$ for some $q\in (1,\infty)$. Then for each smooth path $\lambda(t)$ lying in $U$ with nontangential limit $\lambda(0)=y\in S$ we have
$$
\lim_{t\rightarrow 0}\int_{S'}\tilde{C}'_{p,N}(x',\lambda(t)')n'(x')\psi'(x')d\sigma'(x')=\frac{1}{2}\psi'(y)
$$
$$
+P.V.\int_{S'}\tilde{C}'_{p,N}(x',y')n'(x')\psi'(x')d\sigma'(x')
$$ 
for almost all $y'\in S'$.
\end{theorem}

\ Further as the term $C_{p,N}(x,y)-G(x-y)$ is bounded on $S$ it follows from arguments presented for instance in \cite{MM} that the singular integral
$$
P.V.\int_{S'}\tilde{C}_{p,N}(x',y')n'(x')\psi'(x')d\sigma(x')
$$
defines an $L^{q}$ bounded operator $\Sigma_{S'}:L^{q}(S')\rightarrow L^{q}(S')$.

\ Similarly we have that if $U$, $S$ and $\psi$ are as in Theorem 3 but now $\lambda$ is a smooth path lying in ${\cal{F}}_{p}[N]\backslash(U\cup S)$ with nontangential limit $\lambda(0)=y\in S$ then
$$
\lim_{t\rightarrow 0}\int_{S'}\tilde{C}'_{p,N}(x',\lambda(t)')n'(x')\psi'(x')dm'(x')=-\frac{1}{2}\psi'(y')
$$
$$
+PV\int_{S'}\tilde{C}'_{p,N}(x',y')n'(x')\psi'(x')dm'(x')
$$
for almost all $y\in S'$.

\ It may easily be determined by adapting arguments from \cite{MM} and elsewhere to the situation described here that the operator $\frac{1}{2}I+\Sigma_{S'}$ acting on $L^{q}(S')$ is a projection onto the generalized Hardy space $H^{q,+}$ of left monogenic sections on $U$ whose nontangential maximal function on $S'$ belongs to $L^{q}(S')$. Here $I$ is the identity operator acting on $L^{q}(S')$. However for $p\leq n-2$ the operator $-\frac{1}{2}I+\Sigma_{S'}$ is not necessarily a projection operator acting on $L^{q}(S')$. 

\ To see this recall that in Section 3 the fundamental domain ${\cal{F}}_{p}[N]$ is set up so that for $p\leq n-2$ nontrivial open subsets of $\partial H^{+}(\mathbb{R}^{n})$ belong to $\partial{\cal{F}}_{p}[N]$. Suppose that $w$ belongs to such an open subset then in general $$
\lim_{y\rightarrow w}\int_{S'}\tilde{C}'_{p,N}(x',y')n'(x')\psi'(x')d\sigma(x')
$$
 need not be zero. 

\ To overcome this situation let us instead of considering the fundamental domain ${\cal{F}}_{p}[N]$ let us instead consider the fundamental domain of $\mathbb{R}^{n}\backslash\Gamma_{p}[N]$ which contains ${\cal{F}}_{p}[N]$. We shall denote this fundamental domain by ${\cal{G}}_{p}[N]$. It should be noted that now $w\in{\cal{G}}_{p}[N]$ and that if $y(t)$ is a path in ${\cal{G}}_{p}[N]$ that tends to infinity then $\lim_{t\rightarrow \infty}\int_{S}\tilde{C}_{p,N}(x,y(t))n(x)\psi(x)d\sigma(x)$ is zero. 

\ Let us now introduce the generalized Hardy space $H^{q,-}(S')$ of left monogenic sections defined on ${\cal{N}}_{p}[N]\backslash(U'\cup S')$ whose nontangential maximal function belongs to $L^{q}(S')$. Further ${\cal{N}}_{p}[N]$ is the conformally flat manifold obtained through the factorization $\mathbb{R}^{n}\backslash\Gamma_{p}[N]$. It may now be determined that the operator $-\frac{1}{2}I+\Sigma_{S'}$ is a projection operator from $L^{q}(S')$ onto $H^{q,-}(S')$. Consequently we have:

\begin{theorem}
Suppose that $S$ is as in Theorem 3. Then for $q\in (1,\infty)$ and $p\leq n-2$
$$
L^{q}(S')=H^{q,+}(S')\oplus H^{q,-}(S').
$$
\end{theorem}

\ One can go further than this and set up operators of Calderon-Zygmund type in this context. Suppose that $\phi:\mathbb{R}^{n}\rightarrow Cl_{n}$ is an odd, smooth function that is homogeneous of degree zero. Now consider the kernel $K(x-y)=\frac{\phi(x-y)}{\|x-y\|^{n}}$. For $S$ a strongly Lipschitz surface lying in  $\mathbb{R}^{n}$ then provided $K(x-y)$ satisfies the usual cancellation property described in \cite{MM} and elsewhere the singular integral $P.V.\int_{S}K(x-y)n(x)\psi(x)d\sigma(x)$ defines an operator $K_{S}:L^{q}(S)\rightarrow L^{q}(S)$ of Calderon-Zygmund type. From Lemma 1 and Proposition 3 it now follows that as $K(x)$ is homogeneous of degree $1-n$ the series 
$$
K_{p,N}(x,y):=\sum\limits_{T\in{\cal{T}}_{p}[N]}\sum\limits_{M:\Gamma_{p}[N]\backslash{\cal{T}}_{p}[N]}K(TM<x>-y)J(M,x)
$$ 
is uniformly convergent on each compact subset of $H^{+}(\mathbb{R}^{n})\backslash\cup_{M\in\Gamma_{p}[N]}\{M<y>\}$. Further, by the same arguments used to establish Proposition 5 we may determine the following automorphic invariance of the kernel $K_{p,N}(x,y)$. 

\begin{proposition}
For each $L\in\Gamma_{p}[N]$ we have $J(L,x)K_{p,N}(x,y)=K_{p,N}(L<x>,y)$.
\end{proposition}

\ Consequently we have:

\begin{theorem}
Suppose that $p<n-2$. Suppose also that $S$ is a strongly Lipschitz surface lying in ${\cal{F}}_{p}[N]$. Then the operator $K'_{S'}$ defined by the singular integral 
$P.V.\int_{S'}K'(x',y')n'(x')\psi'(x')d\sigma(x')$ is $L^{q}$ bounded for $q\in(1,\infty)$.
\end{theorem}

\section{$k$-Hypergenic functions}

\ In this section we turn to look at analogous results in the hyperbolic setting. In the first subsection we introduce hypermonogenic Eisenstein series that project to hypermonogenic sections defined on a particular spinor bundle introduced here. However, the results produced in this first subsection automatically carry over for $k$-hypergenic functions and even for $k$-hyperbolic harmonic Eisenstein series. For this reason we treat all cases together in this subsection. In the second subsection we focus on introducing the fundamental solutions on ${\cal{M}}_{p}[N]$ of hypermonogenic sections and hyperbolic harmonic functions together with some of their basic properties.

\subsection{$k$-hypergenic Eisenstein series}

The simplest example of a $k$-hypergenic function is the constant function $F(x)=1$. As pointed out in our preliminary section if $f(y)$ is $k$-hypergenic in the variable $y=M<x>$ then $K_{k}(M,x)f(M<x>)$ is a $k$-hypergenic function in the variable $x$. Consequently upon applying this conformal weight factor to the constant function $F(y)=1$ and Lemma 1 one may now introduce $k$-hypergenic Eisenstein series in upper half space as follows:

\begin{proposition}\cite{ConKraRy}
For $p < n$ and arbitrary real $k$ with $k < n-p-2$ the series
\begin{equation}\label{khypereisen}
\varepsilon_{k,p,N}(x) :=\sum\limits_{M:\Gamma_{p}[N]\backslash{\cal{T}}_{p}[N]}K_{k}(M,x).
\end{equation}
is uniformly convergent on $H^{+}(\mathbb{R}^{n})$ and defines a k-hypergenic function.
\end{proposition}

\ Notice that in all cases where $k < -1$, these series converge even for $p=n-1$.

\par\medskip\par

In complete analogy to the proof of Proposition 6 one can show that
$$
\varepsilon_{k,p,N}(x) =   K_{k}(M,x)\varepsilon_{k,p,N}(M<x>) \quad\quad \forall M \in \Gamma_p[N].
$$
So the series $\varepsilon_{k,p,N}(x)$ defines a $k$-hypergenic Eisenstein series for the group $\Gamma_{p}[N]$.
Furthermore, for $N \ge 3$ we have $\lim\limits_{x_n \rightarrow +\infty} \varepsilon_{k,p,N}(e_n x_n) = 1$. This ensures, that the series $\varepsilon_{k,p,N}(x)$ are non-vanishing functions. 

\par\medskip\par

Following Proposition 1 and \cite{ConKraRy} if   
$$
f(x) = K_{k}(M,x)f(M<x>) \;\forall M \in \Gamma_p[N],
$$
is a $k$-hypergenic automorphic form then the function 
$$
g(x) := \frac{f(x) e_n}{x_n^k}=K_{-k}(M,x)g(M<x>) \;\forall M \in \Gamma_p[N].
$$
is a $-k$-hypergenic automorphic form.
 
This allows us readily to construct non-vanishing $k'$-hypergenic Eisenstein series for positive $k'>1$ from the $k$-hypergenic Eisenstein series of negative $k<-1$, simply by forming
$$
E_{-k,p,N} (x) := \frac{\varepsilon_{k,p,N}(x) e_n}{x_n^k}.
$$
The series $E_{-k,p,N} (x)$ then satisfy the transformation law
$$
E_{-k,p,N} (x) := K_{-k}(M,x)E_{-k,p,N}(M<x>)
$$    
for all $M \in \Gamma_p[N]$. Since the original series $\varepsilon_{k,p,N}(x)$ are non-vanishing functions, the series  $E_{-k,p,N} (x)$ do not vanish, either. 

In particular, this construction provides us with non-trivial hypermonogenic Eisenstein series, which we obtain by putting $k=-n+2$, i.e. 
$$
E_{n-2,p,N} (x) := \frac{\varepsilon_{-n+2,p,N}(x) e_n}{x_{n}^{2-n}} = \sum\limits_{M:\Gamma_p[N]\backslash{\cal{T}}_{p}[N]} \frac{1}{x_n^{2-n}}K_{-n+2}(M,x)e_{n}.
$$
These satisfy in particular 
$$
E_{n-2,p,N} (x) := K_{2-n}(M,x) E_{n-2,p,N}(M<x>)
$$
for all $M \in \Gamma_p[N]$. 

\par\medskip\par

Let us denote by $E_k$ the particular spinor bundle over ${\cal{M}}_{p}[N]$ constructed by
making the identification $(x,X) \leftrightarrow (M<x>,K_{k}(M,x)X)$ for every $M \in \Gamma_p[N]$ where
$x \in H^+(\mathbb{R}^n)$ and $X \in Cl_n$. If an $E_{k}$ valued section defined on a domain $U'$ of ${\cal{M}}_{k}[N]$ lifts to a $k$-hypergenic function on the covering set $U$ of $U'$ then that section is called a left $k$-hypergenic section.

\ We may now state:
  
The projection map applied to $\varepsilon_{k,p,N}$ induces a well-defined non-vanishing $k$-hypergenic
section with values in the spinor bundle $E_k$.  

\ By similar arguments to those used to introduce $k$-hypergenic Eisenstein series in this section one may determine the following:

\begin{theorem}
For any positive integer $N$, for $p<n$ and an arbitrary real $k$ with $k<n-p-2$ the series
$$
\mu_{k,p,N}(x):=\sum\limits_{M:\Gamma_{p}[N]\backslash{\cal{T}}_{p}[N]}L_{k}(M,x)
$$
is uniformly convergent on $H^{+}(\mathbb{R}^{n})$ and defines a $k$-hyperbolic harmonic function satisfying $\mu_{k,p,N}(x)=L_{k}(M,x)\mu_{k,p,N}(M<x>)$.
\end{theorem}

\ Thus for this range of $p$ we have introduced $k$-hyperbolic harmonic Eisenstein series.

\ Let us denote by $B_{k}$ the particular bundle over ${\cal{M}}_{k}[N]$ constructed by making the identification $(x,X)$ with $(M<x>,L_{k}(M,x)X)$ for each $M\in\Gamma_{p}[N]$, with $x\in H^{+}(\mathbb{R}^{n})$ and $X\in Cl_{n}$. We shall call a $B_{k}$ valued section defined on an open subset $U'$ of ${\cal{M}}_{k}[N]$ a $k$-hyperbolic harmonic section if it lifts to a $k$-hyperbolic harmonic function on the lifting of $U$. 

\ It follows from Theorem 4 that the Eisenstein series $\Theta_{k,p,N}(x)$ projects to a well defined $k$-hyberbolic harmonic section defined on $B_{k}$.

\subsection{hypermonogenic and hyperbolic harmonic kernels}

A central aspect in the study of $k$-hypergenic sections on this class of manifolds is again to ask for an
explicit representation of the fundamental solutions on such ${\cal{M}}_{p}[N]$ and for an explicit Cauchy integral formula.

\par\medskip\par

The simplest case is the particular case where $k=n-2$, the case of $(n-2)$-hypermonogenic functions. This is the case that we shall deal with here.

\ Following similar arguments to those used to establish Proposition 3 we have:

\begin{proposition}
For $p=1$ 
\newline
(i) The series 
$$
A_{p,N}(x,y):=\sum\limits_{T_{\underline{m}},T_{-\underline{m}}\in{\cal{T}}_{p}[N]\backslash\{I\}}\sum\limits_{M:\Gamma_{p}[N]\backslash{\cal{T}}_{p}[N]}(K_{n-2}(T_{\underline{m}}M,x)+K_{n-2}(T_{-\underline{m}}M,x))
$$
$$
p(M<x>,y)
$$
is uniformly convergent on any compact subset of $H^{+}(\mathbb{R}^{n})\backslash\cup_{M\in\Gamma_{p}[N]}\{M<y>\}$. Further the kernel $A_{p,N}(x,y)$ satisfies the asymmetry relation $A_{p,N}(x,y)=-\tilde{D}_{p,N}(y,x)$, and $K_{n-2}(L,x)A_{p,N}(L<x>,y)=A_{p,N}(x,y)$ for each $L\in\Gamma_{p}[N]$.
\newline
(ii) The series
$$
B_{p,N}(x,y):=\sum\limits_{T\in{\cal{T}}_{p}[N]}\sum\limits_{M:\Gamma_{p}[N]\backslash{\cal{T}}_{p}[N]}K_{2-n}(TM,x)q(TM<x>,y)
$$
is uniformly convergent on any compact subset of $H^{+}(\mathbb{R}^{n})\backslash\cup_{M\in\Gamma_{p}[N]}\{M<y>\}$. 
\end{proposition}

In order to obtain a Cauchy integral formula for hypermonogenic functions in the context we are considering here let us first note that the kernel $p(x,y)$ is hypermonogenic in the variable $x$.  It follows that for each $M\in\Gamma_{p}[N]\backslash I$ we have that 
$$
P(\int_{\partial V}p(M<x>,y)\tilde{K}_{n-2}(M,x)\frac{n(x)}{x_{n}^{n-2}}f(x)d\sigma(x))=0
$$
 for each function $f$ which is left hypermonogenic in a neighbourhood of the closure of the bounded domain $V$. Further we are assuming that $M<x>\ne y$  It follows that we now have the following version of Cauchy's integral formula.

\begin{theorem}
Suppose that $U$ is a domain in upper half space satisfying $M(U)=U$ for each $M\in \Gamma_{p}[N]$. Suppose also that $f:U\rightarrow Cl_{n}$ is a left hypermonogenic function satisfying $K_{n-2}(M,x)f(M<x>)=f(x)$. Suppose further that $V$ is a bounded subdomain of $U$ and that the closure of $V$ lies in a fundamental domain of $\Gamma_{p}[N]$. Then for each $y\in V$ and for $p<n-2$ 
\begin{equation}
\label{ppseries}
P(f(y))=P\int\limits_{\partial V}\tilde{A}_{p,N}(x,y)\frac{n(x)}{x_{n}^{n-2}}f(x)d\sigma(x).
\end{equation}
\end{theorem} 

To obtain a complete Cauchy integral formula for all of $f(y)$ in this context let us first note from Proposition 1(i) that if $f(x)$ is right hypermonogenic then $e_{n}x_{n}^{n-2}f(x)$ is right hypermonogenic. Further in \cite{sl} it is shown that $q(x,y)$ is right $(2-n)$-hypergenic.  Now consider the integral 
$$
\int_{\partial M^{-1}(V)}e_{n}u_{n}^{n-2}q(u,y)\frac{n(u)}{u_{n}^{n-2}}g(u)d\sigma(u)
$$
 where $u=M<x>$ and $u_{n}$ is the $n$th component of $u$. Further $g$ is a $C^{1}$ function defined on $U$, where $U$ and $V$ are as in Theorem 5. Under a conformal change in variables this integral becomes 
$$
\int_{\partial V}e_{n}q(M<x>,y)\tilde{K}_{2-n}(M,x)n(x)K_{n-2}(M,x)g(M<y>)d\sigma(x).
$$
It follows from Proposition 1(i) that $q(M<x>,y)\tilde{K}_{2-n}(M,x)$ is right $(2-n)$-hypergenic in $x$. Consequently if $U$, $V$, $y$ and $f$ are as in Theorem 5 then 
$$
Q\int_{\partial V}q(M<x>,y)\tilde{K}_{2-n}(M,x)n(x)f(x)d\sigma(x)=0
$$
 for each $M\in \Gamma_{p}[N]\backslash\{I\}$.
Consequently we have 
\begin{theorem}
Suppose $U$, $V$, $y$ and $f$ are as in Theorem 5 and $p<n-2$ then
\begin{equation}
\label{qqseries}
Q(f(y))=Q(\int\limits_{\partial V}\tilde{B}_{p,N}(x,y)n(x)f(x)d\sigma(x))e_{n}.
\end{equation}
\end{theorem}
Combining Theorems 5 and 6 we have
\begin{theorem}
Suppose $U$, $V$, $y$ and $f$ are as in Theorems 5 and 6, and $p<n-2$ then
$$
f(y)=P(\int\limits_{\partial V}\tilde{A}_{p,N}(x,y)\frac{n(x)}{x_{n}^{n-2}}f(x)d\sigma(x))+Q(\int_{\partial V}\tilde{B}_{p,N}(x,y)n(x)f(x)d\sigma(x))e_{n}.
$$
\end{theorem}
 
\ By the same arguments used to prove Proposition 3 we may deduce:

\begin{proposition}
The series
$$
E_{2-n,N}(x,y):=\sum\limits_{T\in{\cal{T}}_{p}[N]}\sum\limits_{M:\Gamma_{p}[N]\backslash{\cal{T}}_{p}[N]}K_{2-n}(TM,x)E(TM<x>,y)
$$
and
$$
F_{2-n,N}(x,y):=\sum\limits_{T\in{\cal{T}}_{p}[N]}\sum\limits_{M:\Gamma_{p}[N]\backslash{\cal{T}}_{p}[N]}\hat{K}_{2-n}(TM,x)F(TM<x>,y)
$$
converge uniformly on each compact subset of $H^{+}(\mathbb{R}^{n})\backslash\cup_{M\in\Gamma_{p}[N]}\{M<y>\}$. 
\end{proposition}

\ By the same arguments used to prove Proposition 5 we can deduce:

\begin{proposition}
The kernels $E_{2-n,N}(x,y)$ and $F_{2-n,N}(x,y)$ satisfy the relationships $K_{2-n}(L,x)E_{k,N}(L<x>,y)=E_{k,L}(x,y)$ and $\hat{K}_{2-n}(L,x)F_{2-n,N}(L<x>,y)=F_{2-n,N}(x,y)$ for each $L\in\Gamma_{k}[N]$.
\end{proposition}

>From results in \cite{sl} one may show that if $\psi\in L^{q}(S)$, with $q\in (1,\infty)$ and $S$ a strongly Lipschitz hypersurface in $H^{+}(\mathbb{R}^{n})$ then the integral 
$$
 y_{n}^{n-2}(\int_{S}E(x,y)n(x)\psi(x)d\sigma(x)-\int_{S}F(x,y)\hat{n}(x)\hat{\psi}(x)d\sigma(x))
$$
defines a hypermonogenic function on $H^{+}(\mathbb{R}^{n})\backslash S$. By the same arguments used in \cite{qber} it may be determined that the integral
$$
y_{n}^{n-2}(\int_{S}E(M<x>,y)\tilde{K}_{2-n}(M,x)n(x)\psi(x)d\sigma(x)
$$
$$
-\int_{S}F(M<x>,y)\hat{\tilde{K}}_{2-n}(M,x)\hat{n}(x)\hat{\psi}(x)d\sigma(x))
$$
defines a hypermonogenic function on $H^{+}(\mathbb{R}^{n})\backslash\partial V$ for each $\psi\in L^{q}(\partial V)$ with $q\in (1,\infty)$ and $S\subset{\cal{F}}_{p}[N]$. Consequently we have:
\begin{proposition}
Suppose $S\subset{\cal{F}}_{p}[N]$ and $\psi\in L^{q}(\partial V)$ with $q\in (1,\infty)$. Further suppose that $p<n-2$. Then the integral 
$$
y_{n}^{n-2}(\int_{S}\tilde{E}_{2-n,N}(x,y)n(x)\psi(x)d\sigma(x)
-\int_{S}\tilde{F}_{2-n,N}(x,y)\hat{n}(x)\hat{\psi}(x)d\sigma(x))
$$
defines a hypermonogenic function on $H^{+}(\mathbb{R}^{n})\backslash\cup_{M\in\Gamma_{p}[N]}\{M<y>\}$. 
\end{proposition}

\ It is straightforward to verify using arguments given in \cite{qber} that the term
$$
K_{n-2}(M,y)(M<y>_{n})^{n-2}(\int_{M<S>}\tilde{E}_{2-n,N}(M<x>,M<y>)n(M<x>)
$$
$$
\psi (M<x>)d\sigma(M<x>)
$$
$$
-\int_{M<S>}\tilde{F}_{2-n,N}(M<x>,M<y>)\hat{n}(M<x>)\hat{\psi}(M<x>)d\sigma(M<x>))
$$
is equal to 
$$
y_{n}^{n-2}(\int_{S}\tilde{E}_{2-n,N}(x,y)n(x)K_{n-2}(M,x)\psi(x)d\sigma(x)
$$
$$
-\int_{S}\tilde{F}_{2-n,N}(x,y)\hat{n}(x)\hat{K}_{n-2}(M,x)\hat{\psi}(x)d\sigma(x)).
$$
for each $M\in\Gamma_{p}[N]$. This establishes a conformal invariance for our Cauchy type integral.  
 
\ Suppose now that $w$ belongs to an open subset of $\partial{\cal{F}}_{p}[N]\cap\partial H^{+}(\mathbb{R}^{n})$ then
$$
\lim_{y\rightarrow w}y_{n}^{n-2}(\int_{S}\tilde{E}_{2-n,N}(x,y)n(x)\psi(x)d\sigma(x)-\int_{S}\tilde{F}_{p,N}(x,y)\hat{n}(x)\hat{\psi}(x)d\sigma(x)=0
$$ 
for each $\psi\in L^{q}(S)$ with $q\in (1,\infty)$.
Further for $y(t)$ an unbounded path in ${\cal{F}}_{p}[N]$ we have
$$
\lim_{t\rightarrow\infty}y_{n}^{n-2}(\int_{S}\tilde{E}_{2-n,N}(x,y(t))n(x)\psi(x)d\sigma(x)-\int_{S}\tilde{F}_{2-n,N}\hat{n}(x)\hat{\psi}(x)d\sigma(x))=0
$$
for each $\psi\in L^{q}(S)$ with $q\in (1,\infty)$.

\ Let us now introduce the Hardy space $H_{n-2}^{q,+}(S')$ of left hypermonogenic sections defined on $U'$ whose nontangential limits lie in $L^{q}(S')$. Further let us introduce the Hardy space $H_{n-2}^{q,-}(S')$ of left hypermonogenic sections defined on ${\cal{M}}_{p}[N]\backslash(U'\cup S')$ whose nontangential limits on $S'$ lie in $L^{q}(S')$. Now by similar arguments to those used in \cite{qber} and in Section 4 we have 

\begin{theorem}
For $q\in (1,\infty)$
$$
L^{q}(S')=H_{n-2}^{q,+}(S')\oplus H_{n-2}^{q,-}(S').
$$
\end{theorem}

\ Following arguments developed to establish Theorem 2 and results in \cite{qber} we also have:

\begin{theorem}
Suppose that $p<n-2$. Further suppose $\psi':{\cal{M}}_{p}[N]\rightarrow E_{2-n}$ is a $C^{1}$ function with compact support. Then for each $y'\in{\cal{M}}_{p}[N]$
$$
\psi'(y')=M_{n-2}y_{n}'^{n-2}(\int_{{\cal{M}}_{p}[N]}\tilde{E}'_{2-n}(x',y')\psi'(x')dm'(x')
$$
$$
-\int_{{\cal{M}}_{p}[N]}\tilde{F}'_{2-n}(x',y')\hat{\psi}'(x')dm'(x').
$$
\end{theorem}

\ We last turn to look at a hyperbolic harmonic kernel. In \cite{sl} it is shown that the function $\frac{1}{\|x-y\|^{n-2}\|\hat{x}-y\|^{n-2}}$ is hyperbolic harmonic in the variable $y$. Consequently the kernel $H(x,y):=\frac{2^{n-2}}{\omega_{n}}\frac{1}{\|x-y\|^{n-2}\|\hat{x}-y\|^{n-2}}$ is  $(2-n)$-hyperbolic harmonic in the variable $y$.

\ Adapting arguments given in \cite{qber} it can be shown that for any $C^{2}$ function $\psi:H^{+}(\mathbb{R}^{n})\rightarrow Cl_{n}$ with compact support one has
$$
\psi(y)=y_{n}^{n-2}\triangle_{2-n}\int_{H^{+}(\mathbb{R}^{n})}H(x,y)\psi(x)dm(x)
$$
for each $y\in H^{+}(\mathbb{R}^{n})$.  

\ Now let us introduce the series 
$$
H_{p,N}(x,y):=\sum\limits_{T\in{\cal{T}}_{p}[N]}\sum\limits_{M:\Gamma_{p}[N]\backslash{\cal{T}}_{p}[N]}L_{2-n}(TM<y>,x)H(TM<y>,x).
$$

\ It follows from Lemma 1 and a straightforward adaptation of the proof of Proposition 3 that this series converges uniformly on any compact subset of $H^{+}(\mathbb{R}^{n})\backslash\cup_{M\in\Gamma_{p}[N]}\{M<y>\}$ for $p<n-3$ and for any $N$. Consequently the kernel $H_{p,N}(x,y)$ is $(2-n)$-hyperbolic harmonic in the variable $y$. Also as each term in the series for $H_{p,N}(x,y)$ is positive it follows that $H_{p,N}(x,y)$ does not vanish for $N=1$ and $N=2$. By similar arguments to those used to establish Proposition 5 it follows that 
$H_{p,N}(y,x)=L_{2-n}(L,y)H_{p,N}(L<y>,x)$ for each $L\in\Gamma_{p}[N]$.

\ Let us now construct a bundle over ${\cal{M}}_{p}[N]$ by making the identification $(x,X)\leftrightarrow (M<x>,\frac{1}{\|cx+d\|^{2n-4}}X)$. We denote this bundle by $B$. We now have in complete analogy to Theorem 10:

\begin{theorem}
Suppose that $p<n-3$, Further suppose $\psi':{\cal{M}}_{p}[N]\rightarrow B$ is a $C^{2}$ section with compact support. Then for each $y'\in{\cal{M}}_{p}[N]$
$$
\psi'(y')=y_{n}'^{n-2}\triangle'_{2-n}\int_{{\cal{M}}_{p}[N]}H'_{p,N}(x',y')\psi'(x')dm'(x')
$$
where $\triangle'_{2-n}$ is the projection to ${\cal{M}}_{p}[N]$ of the operator $\triangle_{2-n}$ and $y'_{n}$ is the projection of $y_{n}$.
\end{theorem}


\begin{thebibliography}{99}

\bibitem{AL} \"O. Akin and H. Leutwiler, On the invariance of the Weinstein equation under M\"obius transformations. {\it Classical and Modern Potential Theory and Applications}, edited by K. GowriSankaran et al, NATO ASI Ser. C,  430, Kluwer, Dordrecht, 1994, 19-29.

\bibitem{a-ord} L. V. Ahlfors, {\it M\"obius transformations in $n$ real variables}. Ordway Lectures Notes, University of Minnesota, 1981.

\bibitem{a} L. V. Ahlfors, M\"{o}bius transformations in $R^{n}$ expressed through $2\times 2$ matrices of Clifford numbers, {\it Complex Variables}, 5, 1986, 215-224.

\bibitem{bds} F. Brackx, R. Delanghe and F. Sommen, {\it{Clifford Analysis}}, Pitman, London, 1982.

\bibitem{bulla} E. Bulla, {\it Fundamental Domains in $\mathbb{R}^k$}, Thesis, University of Arkansas, 2005.

\bibitem{calderbank} D. Calderbank, {\it Dirac operators and Clifford analysis on manifolds with boundary}, Max Planck Institute for Mathematics, Bonn, preprint number 96-131, 1996.

\bibitem{cn} J. Cnops, {\it{An Introduction to Dirac Operators on Manifolds}}, Progress in Mathematical Physics, Birkh\"{a}user, Boston, 2002.

\bibitem{ConKraRy} D. Constales, R. S. Krau{\ss}har, J. Ryan, $k$-hypermonogenic automorphic forms, {\it To appear}

\bibitem{EGM87} J. Elstrodt, F. Grunewald and J. Mennicke,  Vahlen's Group of Clifford matrices and spin-groups, {\it Math. Z.}, 196 1987, 369-390.

\bibitem{EGM88} J. Elstrodt, F. Grunewald and J. Mennicke, Arithmetic applications of the hyperbolic lattice point theorem, {\it Proc. London Math. Soc.} 57, 1988, 239-288.

\bibitem{EGM90} J. Elstrodt, F. Grunewald and J. Mennicke, Kloosterman sums for Clifford algebras and a lower bound for the positive eigenvalues of the Laplacian for congruence subgroups acting on hyperbolic spaces, {\it Invent. Math.}, 101, 1990, 641-668.

\bibitem{EGM98} J. Elstrodt, F. Grunewald and J. Mennicke. {\it Groups Acting on Hyperbolic Space}. Springer, Berlin-Heidelberg, 1998

\bibitem{slberlin} S.-L. Eriksson-Bique, $k$-hypermonogenic functions, {\it Progress in Analysis, Proceedings of the 3rd International ISAAC Congress}, I, edited by H. Begehr et al, World Sci. Publishing, River Edge, New Jersey, 2003, 337-348.

\bibitem{sl} S.-L. Eriksson, Integral formulas for hypermonogenic functions, {\it Bull. Belg. Math. Soc.}, 11, 2004.

\bibitem{slnewnew} S.-L. Eriksson-Bique, M\"obius transformations in several function classes, Univ. Joensuu Dept. Math. Rep. Ser 7, 213-226, 2004.

\bibitem{fk} H. Farkas and I. Kra, {\it Theta Constants, Riemann Surfaces, and the Modular Group }, American Mathematical Society, Rhode Island, 2001.

\bibitem{Ford} L. R. Ford, {\it Automorphic Functions}, McGraw-Hill, New York, 1929.

\bibitem{Freitag} E. Freitag, {\it Hilbert Modular Forms}, Springer Verlag, Berlin, 1990.

\bibitem{FH} E. Freitag and C. F. Hermann.  Some modular varieties of low dimension. {\it Adv. Math.} {\bf 152} No.2 (2000), 203-287.

\bibitem{GR} I. Gradshteyn, and I. M.~Ryzhik, {\it Table of Integrals, Series and Products}, Academic Press, New York, 1980.

\bibitem{GS2} K.~G\"urlebeck and W.~Spr\"ossig. {\it Quaternionic and Clifford Calculus for
Physicists and Engineers}. John Wiley \& Sons, Chichester-New York, 1997.

\bibitem{Hua} L.-K. Hua, {\it Starting with the Unit Circle}, Springer Verlag, Berlin,1981.

\bibitem{KraHabil}  R. S.~Krau{\ss}har, {\it Generalized Analytic Automorphic Forms in Hypercomplex Spaces}, Frontiers in Mathematics, Birkh\"auser, Basel, 2004.

\bibitem{KraBMS} R. S.~Krau{\ss}har, Generalized analytic automorphic forms for some arithmetic congruence subgroups of the Vahlen group on the $n$-dimensional hyperbolic space. {\it Bull. Belg. Math. Soc. Simon Stevin}, 11, 2004, 759--774.

\bibitem{KraRyan1} R. S.~Krau{\ss}har and J. Ryan, Clifford and Harmonic Analysis on Cylinders and Tori, {\it Revista  Matem\'atica Iberoamericana} {\bf 21} (2005), 87--110.

\bibitem{KraRyan2} R. S.~Krau{\ss}har and J. Ryan, Some Conformally Flat Spin Manifolds, Dirac Operators and Automorphic Forms, to appear in {\it Journal of Mathematical Analysis and Applications}. Preprint available at:
http://arXiv.org/abs/math.AP/0212086

\bibitem{KraRyan3} R. S.~Krau{\ss}har, Yuying Qiao and J. Ryan, Harmonic, monogenic and hypermonogenic functions
on some conformally flat manifolds in $\mathbb{R}^n$ arising from special arithmetic groups of the Vahlen group, {\it Contemporary Mathematics}, 370, 2005, 159-173.

\bibitem{Kri88} A. Krieg, Eisenstein series on real, complex and quaternionic half-spaces, {\it Pac. J. Math.}, 133, 1988, 315-354.

\bibitem{Kuiper} N. H. Kuiper, On conformally flat spaces in the large, {\it Ann. Math.}, 50, 1949, 916-924.

\bibitem{LM} H. B. Lawson, M.-L. Michelsohn, {\it Spin geometry}, Princeton University Press, New York, 1989.

\bibitem{Leutw87} H. Leutwiler, Best constants in the Harnack inequality for the Weinstein equation, {\it Aequationes Mathematicae}, 34, 1987, 304--305.

\bibitem{Leutwiler} H. Leutwiler, Modified Clifford analysis, {\it Complex Variables}, 17, 1991, 153-171.

\bibitem{lr} H. Liu and J. Ryan, Clifford analysis techniques for spherical pde, {\it Journal of Fourier Analysis and its Applications}, 8 2002, 535-564.

\bibitem{Maa49} H. Maa{\ss}, Automorphe Funktionen von mehreren Ver\"anderlichen und Dirichletsche Reihen,  {\it Abh. Math. Sem. Univ. Hamb.}, 16, 1949, 53-104.

\bibitem{ma} M. Markel, Regular functions over conformal quaternionic manifolds, {\it Commentationes Mathematicae Universitatis Carolinae}, 22, 1981, 579-583.

\bibitem{MM} M. Mitrea, {\it 
Clifford wavelets, singular integrals, and Hardy spaces}, Lecture Notes in Mathematics, 1575, Springer Verlag, Berlin,  1994.

\bibitem{mitrea} M. Mitrea, Generalized Dirac operators on nonsmooth manifolds and Maxwell's equations, {\it Journal of Fourier Analysis and its Applications}, 7, 2001, 207-256.

\bibitem{qber} Y. Qiao, S. Bernstein, S.-L. Eriksson and J. Ryan, Function Theory for Laplace and Dirac-Hodge Operators in Hyperbolic Space, to appear in {\it Journal d'Analyse Math\'ematique}.

\bibitem{r85} J. Ryan, Conformal Clifford manifolds arising in Clifford analysis, {\it Proc. R. Ir. Acad.}, Sect.,  A 85,  1985, 1-23.

\bibitem{sy} R. Schoen and S.-T. Yau, Conformally flat manifolds, Kleinian groups and scalar curvature, {\it Inventiones Mathematica}, 92, 1988, 47-71.

\bibitem{sijuewu} Sijue Wu, Well-posedness in Sobolev spaces of the full water wave problem in 3D, {\it J. Am. Math. Soc.}, 12, 1999, 445-495.

\bibitem{zoll} G. Z\"oll, {\it Ein Residuenkalk\"ul in der Clifford-Analysis und die M\"obiustransformationen f\"ur euklidische R\"aume}, PhD Thesis, RWTH Aachen, 1987.

\end{thebibliography}
\end{document}